\documentclass[12pt]{amsart}
\usepackage{amsmath,amssymb,amsthm,amsfonts,enumitem,graphicx,comment,xparse,tikz,pgfmath,dsfont,tikz,tikz-cd,verbatim,subcaption,hyperref,mathrsfs,stmaryrd,slashed}
\usepackage[margin=1in]{geometry} 

\usepackage[style=numeric,
sorting=nyt,
date=year,
isbn=false,
doi=false,
url=false,
maxbibnames=99,
backend=biber
]{biblatex}
\AtEveryBibitem{\clearfield{pages}}

\usepackage{tikz}
\usetikzlibrary{decorations.markings}

\renewcommand{\phi}{\varphi}

\setenumerate[0]{wide,label=\textnormal{(\alph*)}}
\newtheorem{theorem}{Theorem}[section]
\newtheorem{corollary}[theorem]{Corollary}
\newtheorem{conjecture}[theorem]{Conjecture}
\newtheorem{lemma}[theorem]{Lemma}
\newtheorem{proposition}[theorem]{Proposition}
\theoremstyle{definition}
\newtheorem{definition}[theorem]{Definition}

\theoremstyle{remark}
\newtheorem{remark}[theorem]{Remark}
\numberwithin{equation}{section}

\newcommand{\alp}{\alpha}
\newcommand{\bt}{\beta}
\newcommand{\gmm}{\gamma}

\newcommand{\dlt}{\delta}
\newcommand{\Dlt}{\Delta}
\newcommand{\eps}{\epsilon}
\newcommand{\veps}{\varepsilon}
\newcommand{\kpp}{\kappa}
\newcommand{\lmb}{\lambda}
\newcommand{\Lmb}{\Lambda}

\newcommand{\tht}{\theta}
\newcommand{\Tht}{\Theta}

\renewcommand{\L}{\mathbb L}

\newcommand{\N}{\mathbb N}

\renewcommand{\P}{\mathbb P}

\newcommand{\R}{\mathbb R}

\newcommand{\T}{\mathbb T}

\newcommand{\Z}{\mathbb Z}

\newcommand{\calA}{\mathcal A}

\newcommand{\calC}{\mathcal C}

\newcommand{\calF}{\mathcal F}
\newcommand{\calG}{\mathcal G}

\newcommand{\calK}{\mathcal K}

\newcommand{\calS}{\mathcal S}

\newcommand{\ul}{\underline}

\newcommand{\tlv}{\tilde v}

\newcommand{\tlx}{\tilde x}

\newcommand{\htP}{{P_{b}}}

\newcommand{\hta}{{a_b}}
\newcommand{\htb}{{b_b}}
\newcommand{\htc}{{\brc_b}}

\newcommand{\htx}{{x_b}}

\newcommand{\httht}{{\theta_b}}
\newcommand{\htTht}{{\Theta_b}}

\newcommand{\htmu}{{\mu_b}}
\newcommand{\htnu}{{\nu_b}}
\newcommand{\htxi}{{\xi_b}}

\newcommand{\httau}{{\tau_b}}

\newcommand{\brC}{\bar C}

\newcommand{\brb}{{\bar{b}_p}}
\newcommand{\brc}{\bar c}

\newcommand{\p}{\partial}

\newcommand{\abs}[1]{|#1|}

\newcommand{\norm}[1]{||#1||}

\newcommand{\pdv}[2][]{\frac{\partial #1}{\partial #2}}

\newcommand{\vel}{\xi}

\newcommand{\ulL}{\ul{L}}
\newcommand{\ulcalC}{\ul{\calC}}

\newcommand{\mra}{{a_r}}
\newcommand{\mrb}{{b_r}}
\newcommand{\mrc}{{\brc_r}}
\newcommand{\mrxi}{{\xi_r}}
\newcommand{\mrtht}{{\tht_r}}
\newcommand{\mrT}{{T_r}}
\newcommand{\mrX}{{X_r}}
\newcommand{\mrL}{{L_r}}
\newcommand{\mrulL}{{\ulL_r}}
\newcommand{\mrmu}{{\mu_r}}
\newcommand{\mrnu}{{\nu_r}}
\newcommand{\mrTht}{{\Tht_r}}

\newcommand{\mrtau}{{\tau_r}}
\newcommand{\mrx}{{x_r}}

\newcommand{\ap}{{a_p}}
\newcommand{\bp}{{b_p}}
\newcommand{\cp}{{\brc_p}}
\newcommand{\xip}{{\xi_p}}

\newcommand{\ep}{\eps}
\title[Instability of Self-similar Implosions]{Instability of Self-Similar Implosions for the Compressible Euler Equations}
\author{Serban Cicortas}
\address{Courant Institute of Mathematical Sciences, 251 Mercer Street, New York, NY 10012}
\address{serbancicortas@nyu.edu}
\author{Igor Rodnianski}
\address{Department of Mathematics, Princeton University, Fine Hall, Washington Road, Princeton, NJ 08544}
\address{irod@math.princeton.edu}
\author{Ely Sandine}
\address{Department of Pure Mathematics and Mathematical Statistics, Centre for Mathematical Sciences, Wilberforce Road, Cambridge CB3 0WB}
\address{ebs37@cam.ac.uk}
\date{}
\begin{document}

\begin{abstract}
We prove that all known regular self-similar implosion solutions of the isentropic compressible Euler equations are unstable to shock formation. The perturbations can be chosen to be small in a class determined by the regularity threshold of the background. The proof relies on an amplification mechanism which occurs along the surface of the backward sonic cone from the implosion point. Motivated by our main results, we formulate an implosion censorship conjecture for the compressible Euler equations, stating that classical implosions are unstable to shock formation.
\end{abstract}

\maketitle

\tableofcontents

\section{Introduction}
\label{sec:introduction}
We consider the Euler equations describing the motion of a compressible fluid on $\R_{t,y}^{1+d}$, $d\geq 2$. We let $\rho(t,y)$ denote the mass density of the fluid, $\xi(t,y)$ the velocity vector field, $p(t,y)$ the pressure and $\veps(t,y)$ the internal energy. The following equations express local conservation of mass, momentum and energy:
\begin{equation}
\label{eq:euler_original}
\begin{cases}
    \p_t \rho+\nabla \cdot (\rho \xi)=0,\\  
    \p_t (\rho \xi)+\nabla \cdot (\rho \xi\otimes \xi)+\nabla p=0,\\
    \p_t(\rho \veps+\frac{1}{2}\rho \xi^2)+\nabla\cdot \left((\rho \veps+\frac{1}{2}\rho \xi^2+p)\xi\right)=0.
\end{cases}
\end{equation}
The system is closed by imposing the ideal gas relation
\begin{equation*}
    \veps=\frac{p}{\rho}\frac{1}{\gmm-1}
\end{equation*}
for an adiabatic index $\gmm>1$. The entropy is given by $S(p \rho^{-\gmm})$ for an increasing function $S$.

The Euler equations are a quasilinear system of conservation laws which possess a well-developed theory of local existence and uniqueness for strong solutions. In particular, their local well-posedness in the family of Sobolev spaces $H^{s}(\R^d)$, $s>\frac{d}{2}+1$, is classical (\cite{katoCauchyProblemQuasilinear1975a}). On the other hand, it is well-known that strong solutions of this type can form singularities.

The first type of singularity, observed already in the 19th century by Challis and Stokes, was shock formation, in which the flow steepens and eventually breaks down via blowing up of its gradient $\norm{\nabla (\rho,\xi,p)}_{L^{\infty}}$. This is a mild singularity, in the sense that one typically expects to be able to develop the solution into a discontinuous shock wave which satisfies \eqref{eq:euler_original} in a weak sense and is unique under appropriate conditions. The theory of conservation laws in the setting of BV solutions is well-developed in the case of $d=1$, see for instance \cite{dafermosHyperbolicConservationLaws2010} and the references therein. However, a corresponding mathematical theory of weak solutions in higher spatial dimensions does not exist yet.

It turns out that in higher dimensions $d>1$ the compressible Euler equations also admit another type of singularities -- implosions. These are stronger, more violent singularities, characterized by infinite concentration and amplification of 
density, velocity and pressure: $\norm{(\rho,\xi,p)}_{L^{\infty}}$.
The construction of implosion singularities arising from regular initial
data is very recent, see \cite{merleImplosionCompressibleFluid2022d,buckmasterSmoothImplodingSolutions2025b,shao2025blowup3dcompressiblenavierstokes,jenssenAmplitudeBlowupRadial2020}.

Beyond the existence of shocks and
implosions there are no general results on classification of singularities to the Euler equations for $d\geq 2$, even in the case of constant entropy and radial symmetry.

The goal of this paper is to show that implosion singularities arising from regular data are unstable. That by itself perhaps is not surprising as previous constructions of implosions indicated a presence of a certain hitherto unknown number of linear unstable directions. 
Here, however, we show that implosions are unstable precisely to shock formation. That is, small 
perturbations of initial data of an implosion solution lead to an earlier shock formation. We identify and exploit an amplification mechanism which occurs along the surface of a backward sonic cone from the implosion point. We then further show that this amplification leads to a shock formation.

These results together with certain analogies between fluid and general relativistic phenomena, specifically between implosions and 
naked singularities, shocks and black hole formation, motivate us to formulate a conjecture arguing that general implosions should always be unstable to shock formation.

The paper will consider 
precisely the class of implosions which are radial, isentropic and self-similar. However, both the amplification and shock formation mechanisms used in this paper hold in the absence of spherical symmetry and the conclusions here can undoubtedly be generalized.

We will now discuss the geometry associated to solutions of the compressible Euler equations, including that of the self-similar 
implosions. 

\subsection{Causal Structure}

Strong solutions to the Euler equations obey finite speed of propagation. We will use the language of Lorentzian geometry to discuss domains of dependence.

The speed of sound is given by 
\begin{equation*}
    c=\sqrt{\frac{\gmm p}{\rho}}.
\end{equation*}
We adopt the convention that the acoustical metric is given by
\begin{align*}
  g&=-dt^2+c^{-2}(dx^i-\vel^i dt)^2.
\end{align*}
The material derivative vector field,
\begin{equation*}
    T:=\p_t+\xi^i\p_{x^i},
\end{equation*}
is a unit-length time-like vector field and defines a time orientation. The acoustical metric and this time orientation describe the causal structure of strong solutions. We note that it is only the conformal class of the acoustic metric which is physically meaningful and normalization conventions vary.

A classical solution is defined on $\calA\subset \R^{1+d}$ which is assumed to be globally hyperbolic with respect to a Cauchy hypersurface $\calK$. A singular point is $p\in \partial \calA$ such that the fluid variables $(\rho,\xi,p)$ do not extend to $p$ in a $C^{1}$ manner. An implosion singularity is $p\in \partial A$ so that $(\rho,\xi,p)$ do not extend to $p$ in a $C^{0}$ manner. Given a singular point $p$, we define its causal past $J^{-}(p)$ as
\begin{equation*}
    J^{-}(p):=\left\{p\right\}\cup \left\{q\in \calA\colon \exists \gmm(t')\colon [0,1)\to\calA,\quad \gmm(0)=q,\quad \lim_{t'\to 1}\gmm(t')=p\right\},
\end{equation*}
where $\gmm$ is assumed to be a $C^{1}$ future-pointing causal curve. We note that this definition does not require extending $g$ to $\partial A$, and hence is applicable to an implosion singularity. The \emph{sonic cone} of a singular point $p$ is the boundary $\partial J^{-}(p)\subset \bar{\calA}$. The sonic cone divides $\calA$ into an interior and exterior region. In particular, given initial data on $\calK$ which develops a singularity at $p$, the singularity is unaffected by exterior perturbations, i.e. those supported on $\calK\setminus J^{-}(p)$.

If one leaves the setting of strong solutions and considers the evolution of shock waves, the above causal considerations are violated. In this case, the Lorentzian metric jumps across the hypersurface of discontinuity. The Rankine-Hugoniot conditions imply that the shock speed is supersonic relative to the fluid ahead of the shock. This is equivalent to the hypersurface being spacelike with respect to the metric ahead of the shock. As such, any sonic cone can be in principle invaded by an exterior shock. Nonetheless, as another consequence of the Rankine-Hugoniot conditions, shocks do not violate Galilean causality. In particular, given a singular point $p$, its sonic cone can only be invaded by shocks forming at $q$ for $t(q)<t(p)$.

The set of all singular points is the \emph{singular boundary}, denoted by $\calS\subset \partial \calA$. From the perspective of classical solutions, every point of $\calS$ is a first singularity. On the other hand, in the setting of weak solutions, the relevant first singular set is its past boundary $\partial_- \calS$, generically a codimension $2$ subset of $\R^{1+d}$. One expects a shock wave to form at this set and propagate into the domain of the classical solution.
\subsection{Self-Similar Implosions}

We now discuss the known examples of implosion singularities. We assume in this section that the flow is isentropic, in which case, the first two equations of \eqref{eq:euler_original} form a closed system for $(\rho,\xi)$.

The resulting system admits the 2-parameter scaling symmetry,
\begin{equation*}
  (\rho(t,y),\vel(t,y))\mapsto \left(\left(\frac{\lmb_2}{\lmb_1}\right)^{\frac{2}{\gmm-1}}\rho\left(\frac{t}{\lmb_1},\frac{y}{\lmb_2}\right),\frac{\lmb_2}{\lmb_1}\vel\left(\frac{t}{\lmb_1},\frac{y}{\lmb_2}\right)\right).
\end{equation*}
For any $\bt>0$, one can consider the $1$-parameter subgroup of scalings
\begin{equation}
  \label{eq:self_sim_scaling}
  (\rho(t,y),\vel(t,y))\mapsto \left(\lmb_1^{\frac{2(\bt-1)}{\gmm-1}}\rho\left(\frac{t}{\lmb_1},\frac{y}{\lmb_1^{\bt}}\right),\lmb_1^{\bt-1}\vel\left(\frac{t}{\lmb_1},\frac{y}{\lmb_1^{\bt}}\right)\right).
\end{equation}
A self-similar solution is one that is invariant with respect to the above scaling for some $0<\bt<1$. 
 If one also assumes radial symmetry, then performing this symmetry reduction leads to a $\bt$-dependent system of ODEs for the profile. The resulting ODE system admits a symmetry, so it can be transformed into an autonomous system. A sufficiently regular global solution of this ODE system leads to a strong solution of \eqref{eq:euler_original} on the domain $\{(t,y)\in \R^{1+d}\colon t<0\}$.

 The resulting solution will feature an implosion at $(t,r)=(0,0)$. The sonic cone of this point will be tangent to the scaling vector field $t\p_t+\bt y^i\p_{y^i}$. This cone leads to a singular sonic point in the ODE system.

 This ODE system is classical, see for instance Guderley (\cite{guderleyStarkeKugeligeUnd1942a}) and Sedov (\cite{sedovSimilarityDimensionalMethods1959}). Early works studied discontinuous shock solutions featuring implosion singularities. The authors of \cite{merleImplosionCompressibleFluid2022d} re-examined the phase portrait and proved the existence of $C^{\infty}$ self-similar implosions. These solutions demonstrate that smooth data can lead to strong singularities for \eqref{eq:euler_original} and were used as profiles to construct blow-up solutions for the compressible Euler, compressible Navier-Stokes and defocusing nonlinear Schr\"odinger equations (\cite{merleImplosionCompressibleFluid2022e,merleBlowEnergySuper2022b}). For the compressible Euler equations, the numerical study \cite{biasiSelfsimilarSolutionsCompressible2021b} indicates that self-similar solutions of the type constructed in \cite{merleImplosionCompressibleFluid2022d} are unstable to shocks which form before the time of implosion. Further constructions of $C^{\infty}$ self-similar profiles were given in \cite{buckmasterSmoothImplodingSolutions2025b} and \cite{shao2025blowup3dcompressiblenavierstokes}. For implosions which are regular but not $C^{\infty}$, see also \cite{jenssenAmplitudeBlowupRadial2020}.

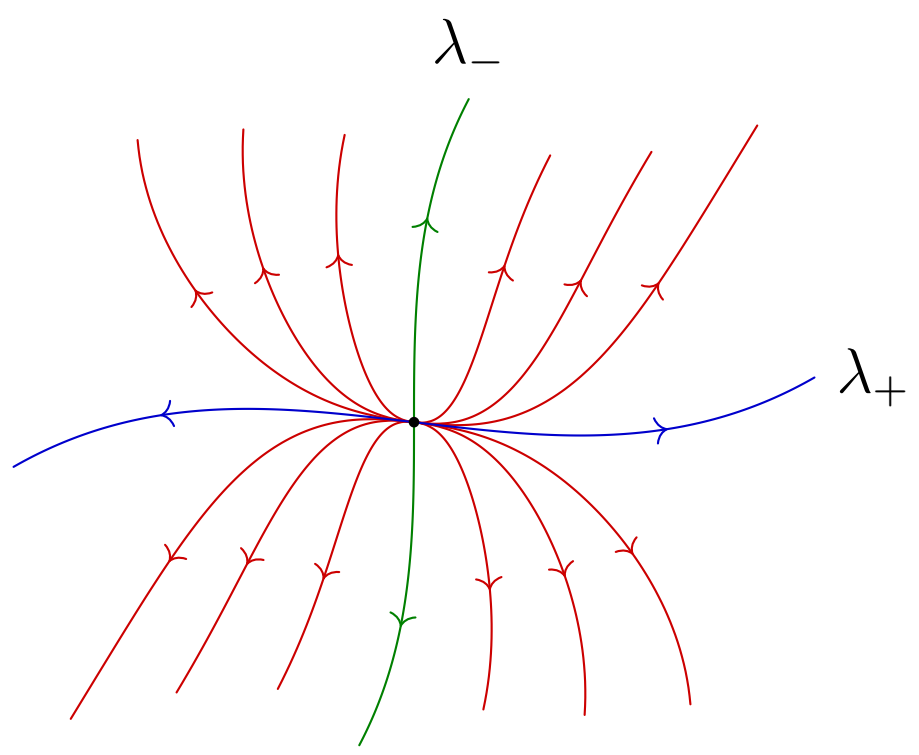
\begin{figure}[htbp]
\centering
\begin{tikzpicture}[
    scale=1,
    line width=0.4pt,
    line cap=round,
    line join=round,
    generic/.style={draw=red!80!black},
    primary/.style={draw=blue!80!black},
    secondary/.style={draw=green!50!black},
    flow/.style={
        postaction={decorate},
        decoration={
            markings,
            mark=at position 0.62 with {\arrow{>}}
        }
    },
    declare function={
        phaseX(\u,\v)=\u+0.035*\v^3;
        phaseY(\u,\v)=\v-0.12*phaseX(\u,\v)
                        +0.033*phaseX(\u,\v)^3;
    }
]

\foreach \sx in {-1,1} {
    \foreach \sy in {-1,1} {
        \foreach \L in {0.68,1.35,2.05} {
            \draw[generic,flow]
                plot[domain=0:1,samples=81,variable=\t]
                ({phaseX(\sx*\L*\t,\sy*1.85*\t^2.5)},
                 {phaseY(\sx*\L*\t,\sy*1.85*\t^2.5)});
        }
    }
}

\foreach \s in {-1,1} {
    \draw[primary,flow]
        plot[domain=0:2.65,samples=81,variable=\t]
        ({phaseX(\s*\t,0)},{phaseY(\s*\t,0)});
}

\foreach \s in {-1,1} {
    \draw[secondary,flow]
        plot[domain=0:2.18,samples=81,variable=\t]
        ({phaseX(0,\s*\t)},{phaseY(0,\s*\t)});
}

\fill (0,0) circle (1pt);

\node[anchor=west]
    at ({phaseX(2.65,0)},{phaseY(2.65,0)})
    {$\lambda_+$};

\node[anchor=south]
    at ({phaseX(0,2.18)},{phaseY(0,2.18)})
    {$\lambda_-$};

\end{tikzpicture}
\caption{Phase portrait of the renormalized ODE system
         near a sonic point.}
\label{fig:phaseportrait_sketch}
\end{figure}

  The ODE system can be renormalized in the neighborhood of a sonic point to be the flow of a vector field near an equilibrium point. For the parameters that we consider, this system will have associated eigenvalues $\lmb_+$ and $\lmb_-,$ satisfying $0<\lmb_+<\lmb_-$. Figure \ref{fig:phaseportrait_sketch} is a qualitative sketch of the phase portrait of the renormalized ODE system near a sonic point.

  The generic integral curve passing through the sonic point, indicated in red, will be tangent to the $\lmb_+$ eigenspace and have regularity $A=\frac{\lmb_-}{\lmb_+}$. We refer to such a solution as a primary direction solution. There is a unique $C^{\infty}$ primary direction solution, indicated in blue. There is a unique integral curve tangent to the $\lmb_-$ eigenspace, indicated in green, which we refer to as the secondary direction solution. The secondary direction solution is $C^{\infty}$. We will consider primary direction solutions. The linearized system about a primary direction solution remembers the generic regularity $A$, which we will refer to as the \emph{regularity threshold}.

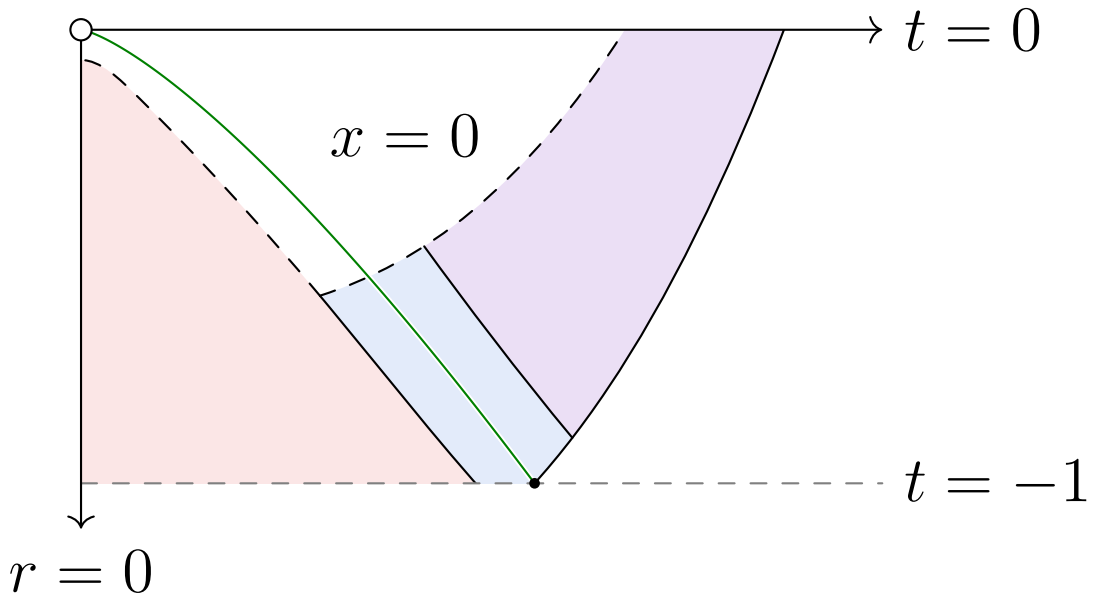
\begin{figure}[htbp]
    \centering
\begin{tikzpicture}[scale=1,line width=0.4pt,
    line cap=round,line join=round]
    \definecolor{figuretwored}{RGB}{225,48,48}
    \definecolor{figuretwoblue}{RGB}{40,103,224}
    \definecolor{figuretwopurple}{RGB}{133,54,194}

    \coordinate (A)  at (0,-4.15);
    \coordinate (E)  at (0,-0.20);
    \coordinate (Raxis) at (0.377977,-0.445788);
    \coordinate (Pm) at (2.80,-3.207790);
    \coordinate (O)  at (3,-3);
    \coordinate (Pp) at (3.25,-2.699161);
    \coordinate (Qm) at (1.580013,-1.759162);
    \coordinate (K)  at (1.633788,-1.823608);
    \coordinate (Qp) at (2.27,-1.432384);
    \coordinate (U)  at (3.596524,0);
    \coordinate (V)  at (4.649233,0);

    \begin{scope}
    \clip (0,-3) rectangle (5.30,0);

    \fill[figuretwored,fill opacity=0.12]
        (A) .. controls (0.933333,-4.131333) and (1.866667,-4.112667) .. (Pm)
        .. controls (2.451239,-2.841027) and (2.000781,-2.265091) .. (K)
        .. controls (1.215184,-1.320039) and (0.796580,-0.860765) .. (Raxis)
        .. controls (0.251984,-0.320888) and (0.125992,-0.20) .. (E)
        -- cycle;

    \fill[figuretwoblue,fill opacity=0.13]
        (Pm) .. controls (2.950000,-3.062363) and (3.100000,-2.894046) .. (Pp)
        .. controls (2.896000,-2.273953) and (2.544000,-1.810116) .. (Qp)
        .. controls (2.040004,-1.580518) and (1.810009,-1.685024) .. (Qm)
        .. controls (1.597938,-1.780563) and (1.615863,-1.802045) .. (K)
        .. controls (2.000781,-2.265091) and (2.451239,-2.841027) .. (Pm)
        -- cycle;

    \fill[figuretwopurple,fill opacity=0.16]
        (Pp) .. controls (3.716411,-2.093183) and (4.182822,-1.230327) .. (V)
        -- (U)
        .. controls (3.154349,-0.701541) and (2.712175,-1.147591) .. (Qp)
        .. controls (2.544000,-1.810116) and (2.896000,-2.273953) .. (Pp)
        -- cycle;

    \draw
        (A) .. controls (1.549744,-4.119005) and (3.099489,-4.088010) .. (V);

    \draw
        (Pm) .. controls (2.451239,-2.841027) and (2.000781,-2.265091) .. (K)
        .. controls (1.615863,-1.802045) and (1.597938,-1.780563) .. (Qm);
    \draw[dashed]
        (Qm) .. controls (1.179335,-1.280787) and (0.778656,-0.842996) .. (Raxis)
        .. controls (0.251984,-0.320888) and (0.125992,-0.20) .. (E);
    \draw (Pp) .. controls (2.896000,-2.273953) and (2.544000,-1.810116) .. (Qp);

    \draw[dashed]
        (Qm) .. controls (2.252184,-1.542490) and (2.924354,-1.066445) .. (U);

    \end{scope}

    \draw[dashed,black!50] (0,-3) -- (5.30,-3)
        node[right,black] {$t=-1$};

    \draw[green!50!black,preaction={draw=white,line width=1.1pt}]
        plot[domain=0:1,samples=161,variable=\s]
        ({3*\s^3},{-3*\s^4});
    \node[anchor=west] at (1.50,-0.70) {$x=0$};

    \draw[->] (0,0) -- (5.30,0) node[right] {$t=0$};
    \draw[->] (0,0) -- (0,-3.30) node[below] {$r=0$};
    \filldraw[fill=white] (0,0) circle (2pt);
    \fill (O) circle (1pt);
\end{tikzpicture}
    \caption{Subregions of the self-similar implosion solution in the $(t,r)$ plane.}
    \label{fig:intro}
\end{figure}

We prove that a class of regular self-similar solutions of the isentropic Euler equations can be perturbed in radial symmetry to form shocks before the time of implosion. The time of implosion for these background solutions is set to be $t=0$. The perturbations we consider will be initially supported near the intersection of $t=-1$ with the background sonic cone, and the shocks will form near this cone. See Figure~\ref{fig:intro} for a sketch of the relevant features of the $(t,r=\abs{y})$ plane.

The curves in Figure~\ref{fig:intro} are all characteristic. The green curve indicates the sonic cone of the self-similar solution, which divides the background spacetime into an interior and exterior region. The blue region is the most general domain that we consider in our shock formation analysis. Notably, this domain is a tubular neighborhood of the sonic cone, far from $r=0$, and not necessarily restricted to the exterior region. By finite speed of propagation, the perturbations we consider vanish in the red shaded region. Finally, the purple region indicates an outer characteristic domain, where in our final result we will prove that the perturbations considered remain regular.
 
We present three theorems which hold for a general class of self-similar implosions, including the ones constructed in \cite{jenssenAmplitudeBlowupRadial2020}, \cite{merleImplosionCompressibleFluid2022d}, \cite{buckmasterSmoothImplodingSolutions2025b} and \cite{shao2025blowup3dcompressiblenavierstokes}, for any dimension $d\geq2$. In the first result, we consider $C^{0,1}$ perturbations. We show by an ODE argument that the one-sided gradient of the solution blows up along the sonic cone. The theorem can be roughly stated as follows.
\begin{theorem}{(Rough version of Theorem \ref{thm:1})}\label{thm:1_intro}
    Consider an admissible self-similar regular implosion background. Consider a radial $C^{0,1}$ initial data perturbation, which satisfies:
    \begin{itemize}
        \item The perturbation vanishes in the interior of the background sonic cone.
        \item The perturbation is $C^{1,1}$ in the exterior of the background sonic cone.
        \item The one-sided derivative of the perturbation at the background sonic cone is non-vanishing with a certain sign.
    \end{itemize}
    It follows that a shock forms in evolution before the background implosion time.
\end{theorem}

In the second result, we consider $C^{N-1,1}$ perturbations where $2\leq N< A$. The theorem is proved using a PDE argument in a subdomain of the exterior and informally stated as follows.
\begin{theorem}{(Rough version of Theorem \ref{thm:2})}\label{thm:2_intro}
    Consider an admissible self-similar regular implosion background and let $A$ denote the regularity threshold. For $2\leq N< A$, consider a radial $C^{N-1,1}$ initial data perturbation, such that:
    \begin{itemize}
        \item The perturbation vanishes in the interior of the background sonic cone.
        \item The perturbation is $C^{N,1}$ in the exterior of the background sonic cone.
        \item For $0\leq k\leq N-1$, the $k$-th one-sided derivative of the perturbation at the background sonic cone vanishes. For $k=N$, this derivative is non-vanishing with a certain sign.
    \end{itemize}
    It follows that a shock forms in evolution before the background implosion time.
\end{theorem} 
We comment that the construction of smooth implosions in \cite{merleImplosionCompressibleFluid2022d} yields backgrounds with $A$ arbitrarily large, so there is no uniform bound on $A$ for the above result. We note that in the first two results, the data can be taken as small as desired, demonstrating instability of the background.  The resulting instability statement is roughly stated as follows.
\begin{corollary}{(Rough version of Corollary \ref{cor:instability_corollary})}
    Consider an admissible self-similar regular implosion background with associated regularity threshold $A>1$. Then for $0\leq k<\lfloor A\rfloor$, the background has at least $\lfloor A\rfloor-k$ unstable directions in the $C^{(k,1)}$ topology.
\end{corollary}

\begin{figure}[htbp]
    \centering
\begin{tikzpicture}[scale=1.15,line width=0.4pt,
    line cap=round,line join=round]
    \definecolor{mghdred}{RGB}{225,48,48}
    \definecolor{mghdblue}{RGB}{40,103,224}
    \definecolor{mghdpurple}{RGB}{133,54,194}

    \coordinate (A)  at (0,-4.15);
    \coordinate (E)  at (0,-0.20);
    \coordinate (Raxis) at (0.377977,-0.445788);
    \coordinate (Pm) at (2.80,-3.207790);
    \coordinate (O)  at (3,-3);
    \coordinate (Pp) at (3.25,-2.699161);
    \coordinate (Qp) at (2.422699,-1.640760);
    \coordinate (U)  at (3.784148,0);
    \coordinate (V)  at (4.649233,0);
    \coordinate (K)  at (1.633788,-1.823608);
    \coordinate (J)  at (1.70,-1.80);
    \coordinate (S)  at (2.00,-1.900316);

    \begin{scope}
    \clip (0,-3) rectangle (5.30,0);

    \fill[mghdred,fill opacity=0.12]
        (A) .. controls (0.933333,-4.131333) and (1.866667,-4.112667) .. (Pm)
        .. controls (2.451239,-2.841027) and (2.000781,-2.265091) .. (K)
        .. controls (1.215184,-1.320039) and (0.796580,-0.860765) .. (Raxis)
        .. controls (0.251984,-0.320888) and (0.125992,-0.20) .. (E)
        -- cycle;

    \fill[mghdblue,fill opacity=0.13]
        (Pm) .. controls (2.950000,-3.062363) and (3.100000,-2.894046) .. (Pp)
        .. controls (2.958774,-2.349355) and (2.668903,-1.973405) .. (Qp)
        .. controls (2.281799,-1.743738) and (2.140900,-1.829240) .. (S)
        .. controls (1.90,-1.855316) and (1.80,-1.82) .. (J)
        .. controls (1.677929,-1.808166) and (1.655859,-1.816032) .. (K)
        .. controls (2.000781,-2.265091) and (2.451239,-2.841027) .. (Pm)
        -- cycle;

    \fill[mghdpurple,fill opacity=0.16]
        (Pp) .. controls (3.716411,-2.093183) and (4.182822,-1.230327) .. (V)
        -- (U)
        .. controls (3.330332,-0.796123) and (2.876515,-1.309085) .. (Qp)
        .. controls (2.668903,-1.973405) and (2.958774,-2.349355) .. (Pp)
        -- cycle;

    \draw
        (A) .. controls (1.549744,-4.119005) and (3.099489,-4.088010) .. (V);

    \draw
        (Pm) .. controls (2.451239,-2.841027) and (2.000781,-2.265091) .. (K);
    \draw[dashed]
        (K) .. controls (1.215184,-1.320039) and (0.796580,-0.860765) .. (Raxis)
        .. controls (0.251984,-0.320888) and (0.125992,-0.20) .. (E);

    \draw (Pp) .. controls (2.958774,-2.349355) and (2.668903,-1.973405) .. (Qp);

    \draw[dashed]
        (K) .. controls (1.655859,-1.816032) and (1.677929,-1.808166) .. (J);
    \draw[thick,dashdotted]
        (J) .. controls (1.80,-1.82) and (1.90,-1.855316) .. (S);
    \draw[dashed]
        (S) .. controls (2.594716,-1.600315) and (3.189432,-1.043301) .. (U);

    \draw[dashed,green!50!black,preaction={draw=white,line width=1.1pt}]
        plot[domain=0:1,samples=161,variable=\s]
        ({3*\s^3},{-3*\s^4});
    \end{scope}

    \draw[dashed,black!50] (0,-3) -- (5.30,-3)
        node[right,black] {$t=-1$};
    \draw[->] (0,0) -- (5.30,0) node[right] {$t=0$};
    \draw[->] (0,0) -- (0,-3.30) node[below] {$r=0$};

    \filldraw[fill=white] (0,0) circle (2pt);
    \filldraw[fill=white] (S) circle (2pt);
    \node[anchor=south,align=center] at (2.05,-0.85) {shock\\formation};
    \draw[->,shorten >=3pt] (2.05,-0.85) -- (S);
    \fill (O) circle (1pt);

\end{tikzpicture}
    \caption{Solution constructed in Theorem~\ref{thm:3_intro}.}
    \label{fig:intro-mghd}
\end{figure}
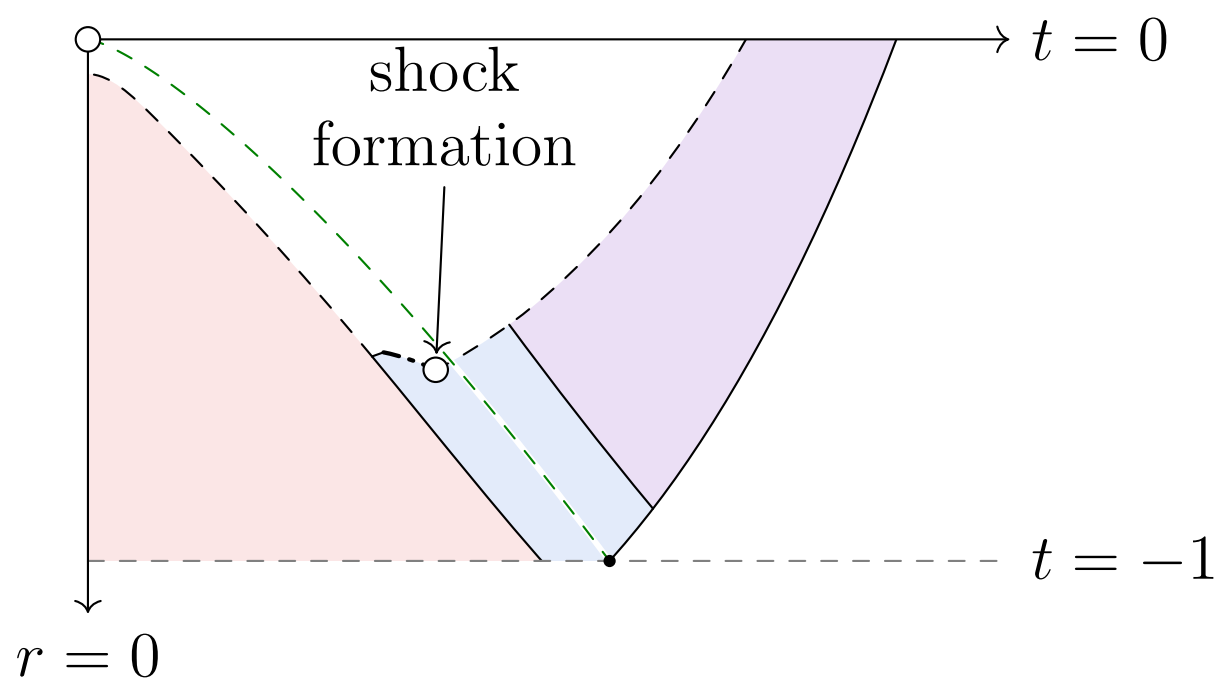

For our third theorem, we consider the case of \textit{smooth} self-similar implosion backgrounds, as the ones constructed in \cite{merleImplosionCompressibleFluid2022d,buckmasterSmoothImplodingSolutions2025b,shao2025blowup3dcompressiblenavierstokes}. We consider $C^{\infty}$ perturbations which quantitatively resemble the $C^{N-1,1}$ perturbations and again lead to a gradient blow-up in evolution. In this result, the data is not restricted to be supported in the exterior region. Under suitable additional assumptions, we prove that the gradient blow-up occurs in the interior of the domain of dependence considered (the blue region in Figure~\ref{fig:intro-mghd}) and so we construct a portion of $\partial_{-}\calS$. Under additional assumptions, the solution remains regular in an outer characteristic region up to $t=0$ (the purple region in Figure~\ref{fig:intro-mghd}).
\begin{theorem}{(Rough version of Theorem \ref{thm:3})}\label{thm:3_intro}
\begin{enumerate}[leftmargin=*]
    \item 
Consider an admissible self-similar smooth implosion background. For $2\leq N< A$, we consider a radial initial data perturbation, which satisfies the following:
\begin{itemize}
    \item The perturbation is $C^{\infty}$ and is supported in a small neighborhood of the background exterior region.
    \item A quantity relating the support of the initial data, the size of its first derivative, and the $C^N$ norm of the initial data is sufficiently small. 
\end{itemize}
It follows that a shock forms in evolution before the background implosion time.
\item Suppose that the following additional hypothesis holds:
\begin{itemize}
  \item The perturbation is supported in a small neighborhood of the background sonic cone.
\end{itemize}
    Then the shock forms, while the solution remains regular on the outer characteristic boundary of the region considered.
  \item Suppose that the following additional hypothesis also holds:
  \begin{itemize}
  \item The background satisfies a suitable exterior stability condition.
\end{itemize}
    Then the solution remains regular in an outer characteristic region up to the background implosion time.
\end{enumerate}
\end{theorem}

As a special case of this result, we may consider a fixed initial data perturbation profile and rescale it to be arbitrarily small in $C^{\lfloor A\rfloor}$, while still achieving shock formation in evolution. In particular, our result shows the instability of a class of self-similar smooth implosions to shock formation, under smooth perturbations which are arbitrarily small in $C^{\lfloor A\rfloor}$. A rough statement is given as follows.
\begin{corollary}{(Rough version of Corollary \ref{cor:smooth_pulse})}
    Consider an admissible self-similar smooth implosion background as in Theorem~\ref{thm:3_intro}. Fix any nontrivial, smooth, compactly supported profile. There exists a sequence of initial data perturbations obtained by rescaling the profile, which converge to $0$ in $C^{\lfloor A\rfloor}$ and satisfy the conditions in Theorem~\ref{thm:3_intro}.
\end{corollary}

\begin{remark}
    [Global solutions] The solutions obtained in Theorem~\ref{thm:3_intro} are defined in the domain of dependence of an outgoing characteristic cone, which extends up to $\{t=0\}.$ We briefly discuss extending the solutions in the exterior of this cone for $\{-1\leq t\leq 0\}.$ In general, for self-similar smooth implosions, constructing a global solution requires physical conditions at infinity at the initial time, such as having a constant state outside of a sufficiently large ball, or sufficiently fast density decay and finite energy. The cut-off procedure needed to enforce these conditions requires an additional notion of exterior stability, see for example \cite{merleImplosionCompressibleFluid2022e, buckmasterSmoothImplodingSolutions2025b}. Under this additional exterior stability condition on the background, the perturbations obtained in Theorem~\ref{thm:3_intro} can also be extended in the exterior of the outgoing characteristic cone for $\{-1\leq t\leq 0\}.$
\end{remark}

For the precise statements of these results, we refer the reader to Section \ref{sec:shock-form-results}.

\begin{figure}[h]
    \centering
\begin{tikzpicture}[
    line width=0.4pt,
    line cap=round,
    line join=round
]
    \definecolor{censorshipblue}{RGB}{40,103,224}
    \definecolor{censorshipred}{RGB}{225,48,48}
    \tikzset{
        boundary/.style={censorshipblue,line width=0.65pt},
        shock/.style={censorshipred,line width=0.65pt},
        horizon/.style={boundary,dashed}
    }

    \foreach \dx/\heading in {
        0/Exterior perturbation,
        6.3/Interior perturbation
    } {
        \begin{scope}[shift={(\dx,0)}]
            \node at (2.05,0.50) {\heading};

            \draw[->,dashed,green!50!black]
                plot[domain=0:1,samples=161,variable=\s]
                ({3.5*\s^3},{-3*\s^4});

            \draw[->] (0,0) -- (4.35,0)
                node[right] {$t=0$};
            \draw[->] (0,0) -- (0,-3.35)
                node[below] {$r=0$};
        \end{scope}
    }

    \begin{scope}
        \coordinate (Sext) at (2.95,-1.70);

        \draw[boundary]
            (0,0) .. controls (0.983333,0)
                         and (1.966667,-1.70) .. (Sext);

        \draw[shock]
            (0,-0.38) .. controls (0.983333,-0.38)
                             and (1.966667,-1.70) .. (Sext);

        \draw[horizon]
            (Sext) .. controls (3.35,-1.70)
                           and (3.90,-1.05) .. (4.10,0);

        \fill[censorshipblue] (Sext) circle (1.2pt);
    \end{scope}

    \begin{scope}[shift={(6.3,0)}]
        \coordinate (Sint) at (2.40,-2.05);

        \draw[boundary]
            (0,-0.50) .. controls (0.80,-0.50)
                             and (1.60,-2.05) .. (Sint);

        \draw[shock]
            (0,-0.85) .. controls (0.80,-0.85)
                             and (1.60,-2.05) .. (Sint);

        \draw[horizon]
            (Sint) .. controls (3.05,-2.05)
                           and (3.70,-1.60) .. (4.10,0);

        \fill[censorshipblue] (Sint) circle (1.2pt);
    \end{scope}

    \foreach \dx in {0,6.3}
        \filldraw[draw=black,fill=white]
            (\dx,0) circle (2pt);
\end{tikzpicture}
    \caption{Plausible evolution of ingoing shocks for solutions in Theorem~\ref{thm:3_intro}.}
    \label{fig:censorship_sketch}
\end{figure}
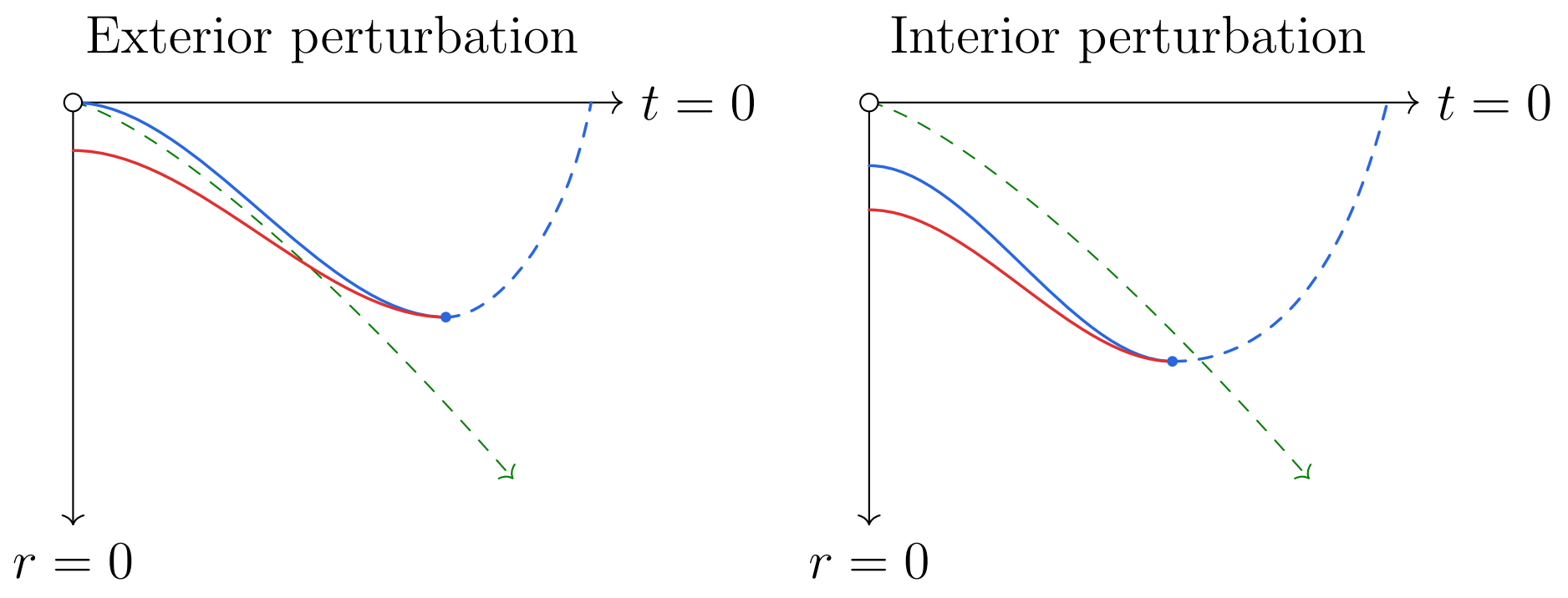

\begin{remark}
    While these results do not address the shock development problem, in Figure~\ref{fig:censorship_sketch} we sketch the plausible evolution of ingoing shocks emanating from the classical solutions in Theorem~\ref{thm:3_intro}. The green curve is the sonic cone of the background. The blue curve represents the boundary of $\calA$ consisting of singular boundary (solid line) and Cauchy horizon (dashed). The red curve indicates the putative shock. The two sketches indicate shocks developing from exterior perturbations and interior perturbations respectively. In the case of an exterior perturbation, the singular boundary is necessarily exterior to the sonic cone. In both cases, one expects the shock to intersect $r=0$ before the background time of implosion, $t=0$.
\end{remark}

\subsection{Discussion of the Proof}
We will work with strong solutions with strictly positive density, which allows us to rewrite the second equation in \eqref{eq:euler_original} as
\begin{equation*}
    \p_t \xi+(\xi\cdot \nabla)\xi+\frac{1}{\rho}\nabla p=0.
\end{equation*}
We restrict to isentropic, radially symmetric solutions, so let $\xi$ denote now the radial velocity and set $r=\abs{y}$. We obtain the following $1+1$ dimensional system:
\begin{equation}
  \label{eq:compressible_euler}
  \begin{cases}
    \p_t\rho+\frac{1}{r^{d-1}}\p_r(r^{d-1}\rho \vel)=0,\\
    \p_t \vel+\vel\p_r\vel+\frac{\p_r p}{\rho}=0.
  \end{cases}
\end{equation}

We begin by reformulating the system using Riemann invariants, which we denote $(a,b)$. We then perform a self-similar renormalization. The new coordinates are $(\tau,x)$ where $\tau$ is a logarithmic time and $x$ is a self-similar coordinate. In this case, $\tau\to\infty$ corresponds to the time of implosion and $x\to-\infty$ corresponds to $r\to0$. The Riemann invariants are replaced by the renormalized variables $(\mra,\mrb)$. We use the convention that $(\hta,\htb)$ correspond to a $\tau$-independent choice of $(\mra,\mrb)$, and represent the background solution. We then decompose the renormalized variables as
\begin{equation*}
  (\mra,\mrb)=(\hta+\ap,\htb+\bp),
\end{equation*}
where $(\ap,\bp)$ now represent the perturbation. To aid the reader, an index of the main coordinates, variables, and definitions introduced in Section~\ref{sec:2} is given in Tables \ref{tab:coordinate_table} and \ref{tab:variable_table}.

\begin{table}[h]
    \centering
    \begin{tabular}{|c|c|}\hline
    $(t,r)$ & Cartesian Radial Coordinates\\\hline
    $(\tau,x)$ & Self-similar Coordinates, \eqref{eq:taux_def}\\\hline
    \end{tabular}
    \caption{Main Coordinates}
    \label{tab:coordinate_table}
\end{table}

\begin{table}[h]
    \centering
    \begin{tabular}{|c|c|}\hline
    $(\rho,\xi)$ & Fluid Variables\\\hline
    $\brc$ & Renormalized Speed of Sound, \eqref{eq:renormalized_speed_sound}\\\hline
    $(a,b)$ & Riemann Variables, \eqref{eq:riemann_variables}\\\hline
    $(\mra,\mrb)$ & Renormalized Riemann Variables, \eqref{eq:renormalized_riemann_variables_def}\\\hline
    $(\hta,\htb)$ & Background Riemann Variables, \eqref{eq:renormalized_riemann_variables}\\\hline
    $(\ap,\bp)$ & Perturbation Variables, \eqref{eq:perturbation_def}\\\hline
    \end{tabular}
    \caption{Main Variables}
    \label{tab:variable_table}
\end{table}

We may assume without loss of generality that $(\hta,\htb)$ has a sonic point that is located at $x=0$. In this case, the linearized system for $(\ap,\bp)$ may be written schematically near $x=0$ as follows, ignoring higher order terms
\begin{equation}
  \label{eq:model_sys}
  \begin{cases}
    (\p_\tau+\kpp x\p_x)\bp-\kpp A \bp+P_{12}\ap=0,\\
    (\p_\tau+2(1-w_0)\p_x)\ap+P_{21}\bp+P_{22}\ap=0,
  \end{cases}
\end{equation}
where $w_0,\kpp,P_{ij}$ are constants depending on the background self-similar profile ($\hta,\htb$) and $A$ is its associated regularity threshold. The background profiles we consider will have $1-w_0>0$, $\kpp>0$ and $A>1$. The condition $1-w_0>0$ corresponds to the speed of sound being positive.

Constants $\kpp$ and $A$ play a crucial role in our analysis. 
They are associated to the sonic cone of the self-similar implosion.
We refer to $\kpp$ as the \emph{surface gravity}, following \cite{merleBlowEnergySuper2022b}. The ``surface gravity'' terminology is used to highlight the connection between the acoustic geometry for this self-similar solution and the geometry of black hole spacetimes, see Remark~3.6 of \cite{merleBlowEnergySuper2022b} and Remark~3.1 of \cite{merleImplosionCompressibleFluid2022e}. In mathematical general relativity, positive surface gravity is well-known to lead to the red-shift of waves, see for instance \cite{dafermosLecturesBlackHoles2013}. As we explain below, the positive surface gravity $\kappa$ will be responsible for the transition between amplification at lower derivative orders and damping at higher derivative orders.

The constant $A$ is the regularity of generic primary direction solutions and also corresponds to the minimum number of derivatives required to achieve coercivity of the linearized operator. This can be seen at the level of the model system \eqref{eq:model_sys}. A first observation is that, if one neglects $\ap$, we have
\begin{equation*}
    \bp|_{x=0}=e^{\kpp A \tau}\cdot \bp|_{(\tau,x)=(0,0)}.
\end{equation*}
In particular, the fluid variables are amplified, analogously to the case of blue-shift effect in general relativity, see already Subsection~\ref{sec:WCC}. Commuting $N$ times with $\p_x$, we see that
\begin{equation*}
  \p_\tau (\p_x^{N}\bp)+\kpp x\p_x (\p_x^{N}\bp)+\kpp(N-A)\p_x^N\bp=0.
\end{equation*}
Restricting to the characteristic $x=0$, we have
\begin{equation}
  \label{eq:growth_heuristic}
  \p_x^N \bp|_{x=0}=e^{\kpp(A-N)\tau}\cdot \p_x^N\bp|_{(\tau,x)=(0,0)}.
\end{equation}
If $N<A$ then the quantity $\p_x^N \bp$ is blue-shifted and hence amplified and if $N>A$ the quantity $\p_x^N \bp$ is red-shifted and hence damped.

In particular, as $A>1$, $\p_x \bp$ will be amplified and the singularities we construct are driven by this linear amplification effect. The non-linear blow-up mechanism is the well-known genuine non-linearity condition of the Euler equations. Theorem~\ref{thm:1_intro} follows from a direct ODE computation showing that, if $\p_x \bp<0$ initially and we consider the characteristic data $\ap|_{x=0}=0$, then $\p_x \bp$ blows up in finite time along $x=0$.

The second and third results concern more regular initial data. For Theorem~\ref{thm:2_intro}, the perturbation data for $x\geq 0$ is given to leading order by
\begin{equation*}
    \bp=b_N x^N+O(x^{N+1}),\qquad b_N<0.
\end{equation*}
We show that if $2\leq N<A$, such data will always lead to a singularity. One way to understand the bound on $N$ is to heuristically use \eqref{eq:growth_heuristic} for $x>0$. In this case, one has
\begin{equation*}
    \bp(\tau,x)\approx e^{A \kpp\tau} b_N x^N+O(x^{N+1}),\qquad \p_x \bp(\tau,x)\approx N e^{(A-1)\kpp\tau}b_N x^{N-1}.
\end{equation*}
From the approximation for $\p_x \bp$, we see that $\p_x \bp=-1$ will occur roughly at a timescale $\tau_*$ which satisfies
\begin{equation*}
    e^{\kpp \tau_*}\approx \frac{1}{\abs{b_N}^{\frac{1}{A-1}}\abs{x}^{\frac{N-1}{A-1}}}.
\end{equation*}
Substituting this timescale into the ansatz for $\bp$ implies that its magnitude is approximately
\begin{equation*}
    \bp\approx \frac{1}{\abs{b_N}^{\frac{1}{A-1}}}x^{\frac{A-N}{A-1}}.
\end{equation*}
Smallness of $\bp$ is necessary to stay within the perturbative regime, so the above equation implies that this can be achieved by picking $x$ small if and only if $N<A$.

In Theorem~\ref{thm:3_intro}, we restrict to smooth background solutions and we consider $C^{\infty}$ perturbations which resemble the perturbations of Theorem~\ref{thm:2_intro}, but are not restricted to being supported in the exterior region. In this case, the above discussion allows the initial data perturbations to be small in $C^N$ if and only if $N<A$.

In the case of our second and third results, one must solve the PDE system and a fundamental difficulty is estimating the solution close to the time of blow-up. To achieve this, consistent with the framework of Christodoulou on multi-dimensional shock formation (\cite{christodoulouFormationShocks3dimensional2007a}), we work in characteristic coordinates and reformulate the problem to that of showing that a certain inverse foliation density of characteristics vanishes. As the system is effectively $1+1$ dimensional, we perform estimates at the level of $L^{\infty}$ norms by integrating along characteristics. For early works on shock formation for systems in planar and radial symmetry, see, for instance, \cite{laxDevelopmentSingularitiesSolutions1964}, \cite{johnFormationSingularitiesOnedimensional1974} and \cite{johnBlowupQuasilinearWave1981a}. We set up the problem in a way similar to the recent works \cite{christodoulouFormationShocksElectromagnetic2016}, \cite{andersonShockFormation1D2025a} and \cite{abbresciaTwoGlobalStability2026}.

We remark that the system \eqref{eq:model_sys} features coupling between $\ap$ and $\bp$. This coupling is controlled by considering the data $\ap=0$ and restricting to a thin characteristic slab. Since the vector field transporting $\ap$ is transverse to the sonic cone, when we integrate in this direction, we gain smallness of order $\abs{x}\ll 1$. The arguments are closed using a hierarchy in the equations to propagate this control. We choose to work in a double-null coordinate system $(u,v)$, straightening both the ingoing and outgoing characteristics. One advantage of these coordinates is that the renormalized characteristic vector fields $\p_u,\p_v$ commute. At the linear level, one can observe that due to the red-shift effect, one will have $x=v e^{\kpp u}$, and so the smallness is translated to $\abs{v}e^{\kpp u}\ll 1$.

As another comment, we do not address the development of shock waves in this work. Nonetheless, we use robust methods and the proofs could in principle be adapted to extract the necessary information to set up this problem. In particular, in the setting of Theorem~\ref{thm:3_intro}, we give a sufficient condition for the shock to form in the interior of the domain of dependence of the data. For sufficiently smooth and generic initial data, the proofs could be adapted to uncover more detailed information, for instance the generic formation of a $C^{\frac{1}{3}}$ cusp. One could also consider the shock development problem, as in the recent works \cite{christodoulouShockDevelopmentSpherical2016} and \cite{buckmasterSimultaneousDevelopmentShocks2022a}. One would expect this shock to travel inwards and censor the background self-similar implosion, see Figure~\ref{fig:censorship_sketch} and Conjecture~\ref{conj}.

The analysis in our work is greatly simplified by working in the setting of constant entropy and radial symmetry. One would expect the instability mechanism to hold in the presence of non-radial perturbations featuring non-trivial entropy and vorticity. We also emphasize that the surface gravity and regularity threshold are geometric and can be defined for non-radial self-similar solutions, in which case one would expect to have similar instability results.

\subsection{Comparison with the Weak Cosmic Censorship Conjecture}\label{sec:WCC}
The geometry of the sonic cone, together with the linear amplification mechanism discussed above, suggests an analogy with the instability of naked singularities to the formation of black holes in general relativity. In the two examples discussed below, naked singularities exhibit a linear instability mechanism (the blue-shift effect), which leads to the formation of trapped surfaces at the nonlinear level. In our setting, the amplification near the sonic cone is responsible for the formation of shocks.

Firstly, in the setting of the spherically symmetry Einstein-Scalar-Field system in $3+1$ dimensions, Christodoulou established in \cite{christodoulouInstabilityNakedSingularities1999a} a suitable formulation of the Weak Cosmic Censorship conjecture, by proving that all naked singularity spacetimes are unstable to black hole formation. The argument in \cite{christodoulouInstabilityNakedSingularities1999a} requires working with a low-regularity class of perturbations. In particular, the $C^{1,\frac{k^2}{1-k^2}}$ examples of $k$-self-similar naked singularities in \cite{christodoulouExamplesNakedSingularity1994} are unstable to $C^{1,\frac{k^2}{16(1-k^2)}}$ perturbations. On the other hand, we remark that for our results the topology used to measure the smallness of the initial data perturbations is limited by the background, but we do not need to restrict to a low-regularity class of perturbations. More precisely, as our third theorem shows, we can prove instability using smooth perturbations that are small in the $C^{N-1,1}$ topology, for $N<A$. 

Secondly, in the setting of the circularly symmetric Einstein-Scalar-Field system in $2+1$ dimensions with a negative cosmological constant, \cite{cicortasWeakCosmicCensorship2026} proved a version of the Weak Cosmic Censorship conjecture using perturbations that are small in a topology of arbitrarily high regularity. In contrast to the present setting, the linear instability mechanism for all the singularities that arise in the model considered by \cite{cicortasWeakCosmicCensorship2026} is so strong, that the regularity of the admissible perturbations is not limited by the background.

Motivated by the analogy between our results and the instability of naked singularities to the formation of black holes, we conjecture that an analogue of the Weak Cosmic Censorship conjecture holds for the compressible Euler equations. We formulate two versions of the conjecture below, for strong solutions and weak solutions. In both cases, we consider the evolution of regular initial data. We call an implosion \textit{classical} if the solution is classical in the causal past of the implosion (in particular, in the case of weak solutions no shocks terminate at the implosion point, which excludes the Guderley solutions \cite{guderleyStarkeKugeligeUnd1942a}, see Figure~\ref{fig:guderley-implosion}). 

\begin{conjecture}[Implosion Censorship Conjecture]\label{conj}
    We consider the compressible Euler equations with variable entropy in $d\geq2$. We conjecture the following:
    \begin{enumerate}
        \item[(1).] (Strong solutions). For generic regular initial data, implosions do not form in classical evolution.
        \item[(2).]  (Weak solutions). For generic regular initial data, weak solutions do not form classical implosions and any implosion of a generic weak solution has a shock that terminates at the implosion.
    \end{enumerate}
    Moreover, for both strong and weak solutions, every classical implosion is unstable to shock formation, in the sense that there exist arbitrarily small perturbations of the initial data whose evolution forms shocks but no classical implosions. Thus, classical implosions are censored by shocks.
\end{conjecture}

We point out that the formulation of the conjecture depends on the initial data regularity and topology considered. It is also reasonable to first consider the conjecture in the case of isentropic spherically symmetric solutions. We further refer the reader to our discussion below, regarding examples of classical implosions outside of irrotationality (\cite{cao-laboraNonradialImplosionCompressible2026}, \cite{chenVorticityBlowupCompressible2025b,chenVorticityBlowup2D2026}) and isentropy (\cite{chenSmoothStableEuler2026b}).

\subsection{Other Implosion Results and their Stability}
Implosions for the compressible Euler and compressible Navier Stokes equations on $\R^{3}$ outside of symmetry and on $\T^{3}$ were constructed in \cite{cao-laboraNonradialImplosionCompressible2026}. For the Euler equations, vorticity blow-ups have been constructed for the isentropic Euler system (\cite{chenVorticityBlowup2D2026}, \cite{chenVorticityBlowupCompressible2025b}). These blow-ups are driven by smooth radial implosions, so we expect that they are unstable to similar perturbations as the ones considered in our work. On the other hand, in the case of the recent examples of non-isentropic self-similar blow-ups in \cite{alexanderSupersonicGravitationalCollapse2025b} and \cite{chenSmoothStableEuler2026b}, the sonic cones of the solutions degenerate, and hence the solutions are not admissible backgrounds for our approach. 

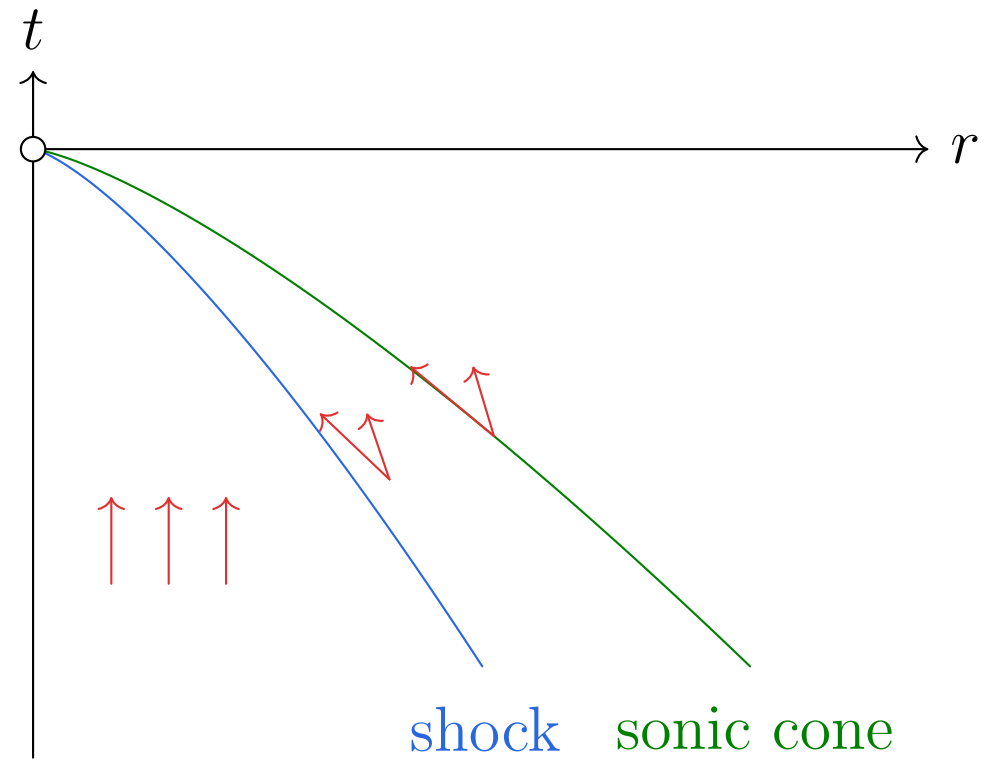
\begin{figure}[htbp]
    \centering
    \begin{tikzpicture}[scale=1.15,line width=0.4pt,
        line cap=round,line join=round]
        \definecolor{guderleyblue}{RGB}{40,103,224}
        \definecolor{guderleyred}{RGB}{225,48,48}

        \draw[guderleyblue]
            plot[domain=0:1,samples=201,variable=\s]
            ({2.60*\s^3},{-3*\s^4});
        \draw[green!50!black]
            plot[domain=0:1,samples=201,variable=\s]
            ({4.15*\s^3},{-3*\s^4});

        \node[guderleyblue,anchor=north] at (2.60,-3.08)
            {shock};
        \node[green!50!black,anchor=north,align=center]
            at (4.15,-3.08) {sonic cone};

        \foreach \r in {0.45,0.78,1.11}
            \draw[->,guderleyred] (\r,-2.50) -- (\r,-2.00);

        \draw[->,guderleyred] (2.05,-1.90) -- (1.65,-1.52);
        \draw[->,guderleyred] (2.05,-1.90) -- (1.92,-1.52);

        \pgfmathsetmacro{\sonicr}{4.15*(1.65/3)^0.75}
        \pgfmathsetmacro{\sonictangent}{0.75*\sonicr/1.65}
        \draw[->,guderleyred]
            (\sonicr,-1.65)
            -- ({\sonicr-0.40*\sonictangent},-1.25);
        \draw[->,guderleyred]
            (\sonicr,-1.65) -- ({\sonicr-0.12},-1.25);

        \draw[->] (0,0) -- (5.15,0) node[right] {$r$};
        \draw[->] (0,-3.50) -- (0,0.45) node[above] {$t$};
        \filldraw[fill=white] (0,0) circle (2pt);
    \end{tikzpicture}
    \caption{Guderley implosion in the $(t,r)$ coordinates.}
    \label{fig:guderley-implosion}
\end{figure}

In addition to regular implosions, there have been many recent investigations into shocks that terminate at implosion points. For recent mathematical works on the existence and formation of the Guderley shock solution, see \cite{jangSelfsimilarConvergingShock2025}, \cite{cialdeaClassicalEulerFlows2025a}.  A sketch of the relative locations of shock front and sonic cone for the Guderley shock solution is provided in Figure \ref{fig:guderley-implosion}. These solutions feature dynamic entropy and hence the local ODE system is different than the one discussed in this paper. The solutions pass through a sonic point and have $A<1$. By the heuristics presented above, this is equivalent to stability near the sonic cone. For a recent result on the stability of these solutions, see \cite{cialdeaNonlinearStabilityGuderley2026}.

We note that for trajectories lying in the Guderley phase portrait, the existing literature suggests that the local stability criterion of $A<1$ at the sonic point is incompatible with the global property of regularly connecting to $r=0$. This is evidence for Conjecture \ref{conj}.

We comment that for the gravitational Euler-Poisson system, smooth self-similar collapses were first constructed numerically, see for instance \cite{larsonNumericalCalculationsDynamics1969}, \cite{penstonDynamicsSelfGravitatingGaseous1969}, \cite{hunterCollapseUnstableIsothermal1977} and \cite{yahilSelfsimilarStellarCollapse1983a}. Such solutions have recently been constructed rigorously, see \cite{guoLarsonPenstonSelfsimilarGravitational2021}, \cite{guoGravitationalCollapsePolytropic2022a}, \cite{sandineHunterSelfsimilarImplosion2024} and \cite{sandineHuntertypeImplosionProfiles2025a}. These self-similar solutions feature a sonic cone. The solutions pass through the sonic point in the secondary direction and one can compute $A<1$. In contrast to the known regular self-similar solutions of Euler, the Larson-Penston solution is nonlinearly stable in radial symmetry, as recently established in \cite{guoNonlinearStabilityLarsonPenston2025a}. We highlight that the local stability of the Larson-Penston solution to shock formation along the sonic cone was first observed by Ori and Piran (\cite{oriSimpleStabilityCriterion1988}).

\subsection{Acknowledgments}
S. Cicortas is partially supported by a Simons Investigator Award (V. Vicol.) I. Rodnianski is supported by a Simons Investigator Award and NSF grant DMS-2405167. E. Sandine is supported by the ERC/UKRI advanced grant SWAT EP/X030644/1. 

\textbf{AI statement.} The authors acknowledge the use of ChatGPT Astra~6 and Claude Opus~5.5 to convert hand-drawn figures to Tikz and to explore numerically the implosion censorship conjecture. S. Cicortas was supported in part by access to OpenAI models provided through the ChatGPT for Academic Researchers program.

\section{Formulations of the System}
\label{sec:2}
We begin this section by defining the Riemann variables and then performing a self-similar renormalization of the polytropic Euler equations. In Subsection~\ref{sec:append-self-simil}, we review properties of the autonomous ODE system that will be relevant for this paper and describe the necessary conditions on the background self-similar profile. In Subsection~\ref{sec:nonlinear-equations}, we introduce the main pair of unknowns, $(a_p,b_p)$, representing the perturbation and their evolution equations. In Subsection~\ref{sec:double-null-form}, we reformulate the system in a double null coordinate system. In Subsection~\ref{sec:append-null-coord}, in preparation for working with the double null formulation, we fix a gauge for the self-similar background and derive asymptotics of the resulting null coordinates. Finally, in Subsection~\ref{sec:pertubation_gauge} we fix a convenient gauge for the perturbation system. 

\subsection{Self-Similar Renormalization}
\label{sec:self-simil-renorm}
The goal of this section is to derive the renormalized system \eqref{eq:renormalized_riemann_variables}
from the system \eqref{eq:compressible_euler}.

We will work with strong solutions, and so if the initial data is isentropic this will be propagated. As such, we will consider the first two equations of \eqref{eq:euler_original} and close the system by assuming the polytropic equation of state
\begin{equation*}
    p=P(\rho),\qquad \rho P'(\rho)=\gmm P(\rho).
\end{equation*}
In this case, the speed of sound is given by
\begin{equation*}
    c=C(\rho),\qquad C(\rho)=\sqrt{P'(\rho)}.
\end{equation*}
It is convenient to work with Riemann variables. We first define the renormalized speed of sound to be
\begin{equation}
\label{eq:renormalized_speed_sound}
  \brc=\brC(\rho),\qquad \brC(\rho):=\int_{0}^{\rho}\frac{C(s)}{s}ds.
\end{equation}
We then define the Riemann variables to be
\begin{equation}
  \label{eq:riemann_variables}
  a=\xi+\brc,\qquad b=\xi-\brc.
\end{equation}
We will always consider $\brc>0$, which leads to the physicality condition
\begin{equation}
  \label{eq:physicality}
  a>b.
\end{equation}
We adopt the following notation for the material derivative and normalized spatial derivative:
\begin{equation*}
  T:=\p_t+\vel\p_r,\qquad X:=c \p_r.
\end{equation*}
We comment that these vector fields are orthogonal with respect to the acoustical metric. We define a scalar $\tht$ as
\begin{equation}
  \label{eq:TX_commutator}    
  [T,X]=\tht\cdot X,\qquad \tht:=\frac{T c-X\xi}{c}.
\end{equation}
We remark that the scalar $\tht$ is at the level of the derivative of the fluid variables and typically non-vanishing.

We define the acoustically null vector fields to be
\begin{equation*}
  \ulL=T-X,\qquad L=T+X.
\end{equation*}
One can then compute that $a,b$ satisfy the following first order system:
\begin{equation}
  \label{eq:first_order_formulation}
  \begin{cases}
    \ulL b=\frac{(d-1)c}{2r}(a+b),\\
    L a=-\frac{(d-1)c}{2r}(a+b).
  \end{cases}
\end{equation}

As the pressure law is polytropic, $c$ will be a linear function of $a,b$. We also let
\begin{equation}
  \label{eq:alp_def}
  \alp:=\frac{\gmm-1}{2},
\end{equation}
so that
\begin{equation*}
  c=\alp \brc.
\end{equation*}
We note that \eqref{eq:first_order_formulation} can be written as
\begin{equation}
  \label{eq:r_renormalized_formulation}
  \begin{cases}
    \ulL(r^{-1}b)+r^{-2}(\xi-\alp\brc)b=r^{-2}\frac{(d-1)\alp}{2}\brc(a+b),\\
    L(r^{-1}a)+r^{-2}(\xi+\alp\brc)a=-r^{-2}\frac{(d-1)\alp}{2}\brc(a+b).
  \end{cases}
\end{equation}
In particular, changing to variables $r^{-1}a,r^{-1}b,r^{-1}\brc,r^{-1}\xi$ removes $r$ from the equation. 

We also note that we can decompose
\begin{equation}
    \label{eq:TX_commutator_renormalized}
    \tht=\tht_1+\tht_2,\qquad \tht_1:=\frac{L(r^{-1}(\alp \brc-\xi))}{2\alp r^{-1}\brc},\qquad \tht_2:=\frac{\ulL(r^{-1}(\alp \brc-\xi))}{2\alp r^{-1}\brc}.
\end{equation}

The next step is to convert the equations to renormalized variables and a self-similar coordinate system. We fix a scaling parameter $\bt>0$. We introduce a self-similar coordinate system $(\tau,x)$ as
\begin{equation}
  \label{eq:taux_def}
  \tau=-\bt \ln(-t),\qquad x=\ln(\frac{r}{(-t)^{\bt}}).
\end{equation}
We define renormalized variables as
\begin{equation}
\label{eq:renormalized_riemann_variables_def}
  a=r(-t)^{-1}\bt \mra,\qquad b=r(-t)^{-1}\bt \mrb.
\end{equation}
We note that we have
\begin{equation*}
  \brc=r(-t)^{-1}\bt\frac{1}{2}(\mra-\mrb),\qquad \xi=r(-t)^{-1}\bt\frac{1}{2}(\mra+\mrb).
\end{equation*}
We note that the $(t,r)$ coordinate vector fields transform as
\begin{equation*}
  (-t)\p_t=\bt(\p_\tau+\p_x),\qquad r\p_r=\p_x.
\end{equation*}
We define the renormalized vector fields
\begin{equation*}
  \mrT=\frac{1}{\bt}(-t)T,\qquad \mrX=\frac{1}{\bt}(-t)X.
\end{equation*}
We also set
\begin{equation*}
  \mrxi=1+\frac{\mra+\mrb}{2},\qquad \mrc=\frac{\mra-\mrb}{2},
\end{equation*}
and we compute that
\begin{equation*}
  \mrT=\p_{\tau}+\mrxi \p_x,\qquad \mrX=\alp\mrc \p_x.
\end{equation*}
We similarly define the renormalized vector fields
\begin{equation*}
  \mrulL=\mrT-\mrX,\qquad \mrL=\mrT+\mrX,
\end{equation*}
and we compute with respect to $(\tau,x)$ coordinates
\begin{equation*}
  \mrulL=\p_{\tau}+(\mrxi-\alp \mrc)\p_x,\qquad \mrL=\p_{\tau}+(\mrxi+\alp \mrc)\p_x.
\end{equation*}
Similarly to \eqref{eq:TX_commutator_renormalized}, we have the identity
\begin{equation}
  \label{eq:renormalized_commutator}
  [\mrT,\mrX]=\mrtht \mrX,\qquad \mrtht:=\mrtht_1+\mrtht_2,\qquad \mrtht_1:=\frac{\mrL(\alp\mrc-\mrxi)}{2\alp\mrc},\qquad \mrtht_2:=\frac{\mrulL(\alp\mrc+\mrxi)}{2\alp\mrc}.
\end{equation}
We note that we can write $\mrulL$ and $\mrL$ in terms of $\mra,\mrb$ as follows
\begin{equation*}
  \mrulL=\p_{\tau}+\left(1-\frac{\alp-1}{2}\mra+\frac{\alp+1}{2}\mrb\right)\p_x,\qquad \mrL=\p_{\tau}+\left(1+\frac{\alp+1}{2}\mra-\frac{\alp-1}{2}\mrb\right)\p_x.
\end{equation*}

In terms of the new coordinates and variables, the system \eqref{eq:r_renormalized_formulation} becomes
\begin{equation}
  \label{eq:renormalized_riemann_variables}  
  \begin{cases}
    \mrulL\mrb+\left(\frac{1}{\bt}-1+\mrxi-\alp\mrc\right)\mrb=\frac{(d-1)\alp}{2}\mrc(\mra+\mrb),\\
    \mrL\mra+\left(\frac{1}{\bt}-1+\mrxi+\alp\mrc\right)\mra=-\frac{(d-1)\alp}{2}\mrc(\mra+\mrb).
  \end{cases}
\end{equation}
We will use the notation convention that $(\hta,\htb)$ represents a solution of \eqref{eq:renormalized_riemann_variables} which is independent of $\tau$. Such a solution corresponds to a self-similar solution of \eqref{eq:compressible_euler}. In this case, the system \eqref{eq:renormalized_riemann_variables} reduces to a system of ODEs. We review relevant features of this ODE system in the next subsection.
\begin{remark}
In the notation of \cite{merleImplosionCompressibleFluid2022d}, the parameters $\alp,\bt$ correspond to $l^{-1},r^{-1}$.
\end{remark}

\subsection{Local theory of the ODE system}
\label{sec:append-self-simil}
In this section, we briefly describe some properties of the autonomous ODE system satisfied by self-similar solution of \eqref{eq:compressible_euler}. We also derive in Lemma~\ref{lem:background_construction} necessary conditions on the background self-similar profile, which in particular are satisfied by the examples of \cite{merleImplosionCompressibleFluid2022d}, according to Proposition~\ref{prop:mrrs}.

We consider $\tau$-independent solutions to \eqref{eq:renormalized_riemann_variables} and obtain the following ODE system:
\begin{equation}
  \label{eq:autonomous_ode_sys}
  \begin{cases}
    \big(1-\frac{\alp-1}{2}\hta+\frac{\alp+1}{2}\htb\big)\p_x\htb+\big(\frac{1}{\bt}-\frac{\alp-1}{2}\hta+\frac{\alp+1}{2}\htb\big)\htb=\frac{(d-1)\alp}{4}\big(\hta^2-\htb^2\big),\\
    \big(1+\frac{\alp+1}{2}\hta-\frac{\alp-1}{2}\htb\big)\p_x\hta+\big(\frac{1}{\bt}+\frac{\alp+1}{2}\hta-\frac{\alp-1}{2}\htb\big)\hta=\frac{(d-1)\alp}{4}\big(\htb^2-\hta^2\big).
  \end{cases}
\end{equation}
The sonic line is the subset of phase space given by
\begin{equation}
  \label{eq:sonic_set}
  1-\frac{\alp-1}{2}\hta+\frac{\alp+1}{2}\htb=0.
\end{equation}
We will now consider a solution which crosses the sonic line. We assume without loss of generality that this crossing occurs at $x=0$.

We now assume $C^{1}$ regularity and expand the solution at $x=0$,
\begin{equation*}
  \hta=a_0+a_1x+o(x),\qquad \htb=b_0+b_1x+o(x).
\end{equation*}
The sonic line condition and the first equation in \eqref{eq:autonomous_ode_sys}, give the system
\begin{equation}
  \label{eq:zeroth_order_sys}
  \begin{cases}
    1-\frac{\alp-1}{2}a_0+\frac{\alp+1}{2}b_0=0,\\
    \big(\frac{1}{\bt}-1\big)b_0=\frac{d-1}{4}\alp\big(a_0^2-b_0^2\big).
  \end{cases}
\end{equation}
We introduce the notation
\begin{equation}
  \label{eq:w0def}
  w_0=-\frac{1}{2}\big(a_0+b_0\big),
\end{equation}
and have from the first equation of \eqref{eq:zeroth_order_sys}
\begin{equation}
  \label{eq:abw_relation}
  a_0=\frac{1-(\alp+1)w_0}{\alp},\qquad b_0=\frac{-1-(\alp-1)w_0}{\alp}.
\end{equation}
The second equation of \eqref{eq:zeroth_order_sys} becomes the following quadratic equation in $w_0$:
\begin{equation}
  \label{eq:w_quadratic}
  (d-1)w_0^2+\left(\frac{1}{\alp}\left(\frac{1}{\bt}-1\right)(\alp-1)-(d-1)\right)w_0+\frac{1}{\alp}\left(\frac{1}{\bt}-1\right)=0.
\end{equation}

We denote the discriminant in \eqref{eq:w_quadratic}
\begin{equation*}
  J_1=\left(\frac{1}{\alp}\left(\frac{1}{\bt}-1\right)(\alp-1)-(d-1)\right)^2-4(d-1)\frac{1}{\alp}\left(\frac{1}{\bt}-1\right).
\end{equation*}
If $J_1\geq 0$, then \eqref{eq:w_quadratic} admits the solutions
\begin{equation*}
  w_{0}=\frac{1}{2}-\frac{\left(\frac{1}{\bt}-1\right)(\alp-1)}{2 \alp(d-1)}\pm \frac{\sqrt{J_1}}{2(d-1)}.
\end{equation*}
The physicality condition \eqref{eq:physicality} can then be expressed as
\begin{equation}
  \label{eq:w0-physicality}
  w_0<1.
\end{equation}

We suppose that $w_0$ is a solution to \eqref{eq:w_quadratic} obeying \eqref{eq:w0-physicality}. We now consider the coefficients $a_1,b_1$. Substituting $x=0$ into the second equation of \eqref{eq:autonomous_ode_sys} yields
\begin{equation*}
  \left(1+\frac{\alp+1}{2}a_0-\frac{\alp-1}{2}b_0\right)a_1+\left(\frac{1}{\bt}+\frac{\alp+1}{2}a_0-\frac{\alp-1}{2}b_0\right)a_0=\frac{(d-1)\alp}{4}\left(b_0^2-a_0^2\right).
\end{equation*}
This can be simplified using \eqref{eq:zeroth_order_sys} to give
\begin{equation}
  \label{eq:a1}
  a_1=-a_0-\left(\frac{1}{\bt}-1\right)\frac{a_0+b_0}{\alp(a_0-b_0)}.
\end{equation}
In terms of $w_0$, we have
\begin{equation*}
  a_1+a_0=\frac{w_0}{1-w_0}\left(\frac{1}{\bt}-1\right).
\end{equation*}
If we expand the first equation of \eqref{eq:autonomous_ode_sys} to $o(x)$, we obtain
\begin{equation*}
  \left(-\frac{\alp-1}{2}a_1+\frac{\alp+1}{2}b_1\right)(b_1+b_0)+\left(\frac{1}{\bt}-\frac{\alp-1}{2}a_0+\frac{\alp+1}{2}b_0\right)b_1=\frac{(d-1)\alp}{2}(a_0a_1-b_0b_1).
\end{equation*}
This equation can be simplified to the following quadratic for $b_1+b_0$:
\begin{equation}
  \label{eq:b1_quad}
  \frac{\alp+1}{2}\left(b_1+b_0\right)^2+(D_1+D_2)(b_1+b_0)+D_3=0,
\end{equation}
where we introduce the notation
\begin{equation*}
  D_1:=\frac{1}{\bt}-1+\frac{(d-1)\alp}{2}b_0,\qquad D_2:=1-\frac{w_0}{1-w_0}\frac{\alp-1}{2}\left(\frac{1}{\bt}-1\right),
\end{equation*}
\begin{equation*}
  D_3:=-\frac{d-1}{2}\alp\left(\frac{1}{\bt}-1\right)a_0\frac{w_0}{1-w_0}+\left(\frac{1}{\bt}-1\right)b_0.
\end{equation*}
We denote the discriminant of the quadratic equation \eqref{eq:b1_quad}
\begin{equation*}
  J_2=\left(\frac{D_1+D_2}{2}\right)^2-\frac{\alp+1}{2}D_3.
\end{equation*}
If $J_2> 0$, we may solve \eqref{eq:b1_quad} and obtain
\begin{equation}
  \label{eq:b1_sols}
  \frac{\alp+1}{2}(b_0+b_1)=-\frac{D_1+D_2}{2}\pm \sqrt{J_2},\qquad 
\end{equation}
For a given solution, we define the surface gravity to be
\begin{equation*}
  \kpp:=-\frac{\alp-1}{2}a_1+\frac{\alp+1}{2}b_1.
\end{equation*}
Provided that $\kpp\neq 0$, we define the regularity threshold to be

\begin{equation*}
  A:=-\kappa^{-1}\left(\frac{\alp+1}{2}(b_0+b_1)+\left(\frac{1}{\bt}-1\right)+\frac{(d-1)\alp}{2}b_0\right).
\end{equation*}
We note that we may express $\kappa$ and $A$ as follows
\begin{equation*}
  \kpp=\frac{\alp+1}{2}(b_0+b_1)+D_2,\qquad A=-\frac{\frac{\alp+1}{2}(b_0+b_1)+D_1}{\frac{\alp+1}{2}(b_0+b_1)+D_2}.
\end{equation*}
In particular, for the two solutions given by \eqref{eq:b1_sols}, we have
\begin{equation}
  \label{eq:kpp_A_formulas}
  \kpp=\frac{D_2-D_1}{2}\pm \sqrt{J_2},\qquad A=\frac{D_2-D_1\mp 2\sqrt{J_2}}{D_2-D_1\pm 2\sqrt{J_2}}.
\end{equation}

\begin{remark}
The formulas for $\kpp,A$ can be derived by linearizing \eqref{eq:autonomous_ode_sys} about a solution $\hta, \htb$. This gives the following linear ODE for $(a,b)$:
\begin{equation}
  \label{eq:linearized_autonomous_system}
  \begin{pmatrix}
    \left(1-\frac{\alp-1}{2}\hta+\frac{\alp+1}{2}\htb\right) & 0\\ 0 & \left(1+\frac{\alp+1}{2}\hta-\frac{\alp-1}{2}\htb\right)
  \end{pmatrix}\p_x
  \begin{pmatrix}
    b\\a
  \end{pmatrix}+\begin{pmatrix}
    P_{11}&P_{12}\\P_{21} &P_{22}
  \end{pmatrix}
  \begin{pmatrix}
    b\\a
  \end{pmatrix}=0,  
\end{equation}
where we introduce the notation
\begin{equation}
  \label{eq:P_entries}
  \begin{cases}
    P_{11}:=\frac{\alp+1}{2}\big(\p_x\htb+\htb\big)+\big(\frac{1}{\bt}-\frac{\alp-1}{2}\hta+\frac{\alp+1}{2}\htb\big)+\frac{(d-1)\alp}{2}\htb,\\
    P_{12}:=-\frac{\alp-1}{2}\big(\p_x\htb+\htb\big)-\frac{(d-1)\alp}{2}\hta,\\
    P_{21}:=-\frac{\alp-1}{2}\big(\p_x\hta+\hta\big)-\frac{(d-1)\alp}{2}\htb,\\
    P_{22}:=\frac{\alp+1}{2}\big(\p_x\hta+\hta\big)+\big(\frac{1}{\bt}+\frac{\alp+1}{2}\hta-\frac{\alp-1}{2}\htb\big)+\frac{(d-1)\alp}{2}\hta.
  \end{cases}
\end{equation}
We note that the coefficient of $\p_x b$ satisfies
\begin{equation}
  \label{eq:kpp_property_a}
  1-\frac{\alp-1}{2}\hta+\frac{\alp+1}{2}\htb=\kpp x+O(x^2),
\end{equation}
and the coefficient $P_{11}$ satisfies
\begin{equation}
  \label{eq:A_property_a}
  P_{11}=-\kpp A+O(x).
\end{equation}
In particular, if we perform a Frobenius expansion for \eqref{eq:linearized_autonomous_system} at $x=0$, we see that the characteristic exponents are $\{0,A\}$, so a generic homogeneous solution will have regularity $x^A$ as $x\to 0$. The Frobenius exponent $0$ is due to the $x$-translation symmetry. Thus, one can view $A$ as the generic regularity of solutions to the ODE system \eqref{eq:autonomous_ode_sys} at the sonic line.
\end{remark}

Returning to \eqref{eq:kpp_A_formulas}, we will always assume the inequality
\begin{equation}
  \label{eq:J2_ineq}
  \sqrt{J_2}<\frac{D_2-D_1}{2}.
\end{equation}
In this case, regardless of choice of sign, we have $\kpp>0$. We comment that it is clear that switching the root of $b_1$ corresponds to $A\to A^{-1}$. In addition, we will always consider the choice of ``$-$'' sign for \eqref{eq:b1_sols} which corresponds to $A>1$.

For a more geometric perspective on the ODE solutions, starting with the ODE system \eqref{eq:autonomous_ode_sys} one can introduce the solution-dependent renormalized coordinate $x'$ on $x>0$ by the equation
\begin{equation*}
    dx'=\frac{1}{1-\frac{\alp-1}{2}\hta+\frac{\alp+1}{2}\htb}dx.
\end{equation*}
One then obtains the renormalized system
\begin{equation}
    \label{eq:renormalized_ode_system}
  \begin{cases}
    \p_{x'}\htb=-\left(\frac{1}{\bt}-\frac{\alp-1}{2}\hta+\frac{\alp+1}{2}\htb\right)\htb+\frac{\left(d-1\right)\alp}{4}\left(\hta^2-\htb^2\right),\\
    \p_{x'}\hta=-\frac{\left(1-\frac{\alp-1}{2}\hta+\frac{\alp+1}{2}\htb\right)}{\left(1+\frac{\alp+1}{2}\hta-\frac{\alp-1}{2}\htb\right)}\left(\frac{1}{\bt}+\frac{\alp+1}{2}\hta-\frac{\alp-1}{2}\htb\right)\hta+\frac{\left(1-\frac{\alp-1}{2}\hta+\frac{\alp+1}{2}\htb\right)}{\left(1+\frac{\alp+1}{2}\hta-\frac{\alp-1}{2}\htb\right)}\frac{\left(d-1\right)\alp}{4}\left(\htb^2-\hta^2\right).
  \end{cases}
\end{equation}
The right-hand side can be viewed as a smooth vector field on the $(\hta,\htb)$ plane, which vanishes at the sonic point. The eigenvalues of the derivative of this vector field can be computed to be (note that our convention differs by a sign compared to \cite{merleImplosionCompressibleFluid2022d})
\begin{equation*}
    \lmb_{\pm}:=\frac{D_2-D_1}{2}\mp \sqrt{J_2}.
\end{equation*}
Under the assumption \eqref{eq:J2_ineq}, we have $0<\lmb_+<\lmb_-$. The integral curves of this vector field passing through the stationary point are sketched in Figure \ref{fig:phaseportrait_sketch}, see also now the accompanying discussion.
The primary and secondary directions correspond to the $-$ and $+$ choice of sign in \eqref{eq:b1_sols} respectively. For a primary direction solution, $A=\lmb_-/\lmb_+$ and for a secondary direction solution, $A=\lmb_+/\lmb_-$.

Returning to the system \eqref{eq:autonomous_ode_sys}, we note that once $b_1$ is determined, one can solve for the higher order Taylor coefficients. Formally, the coefficient of $b_N$ for $N\geq 2$ is given by $\kpp(N-A)$. Thus, the higher Taylor coefficients are uniquely determined provided that $A\notin \Z$. One can assume the solution is $C^{\infty}$ and formally solve for the Taylor coefficients. As the coefficients of the singular ODE are analytic, it follows that the Taylor series constructed in this manner converges absolutely on a neighborhood of $x=0$. This leads to the $C^{\infty}$ primary direction solution defined on a neighborhood of $x=0$.

To construct the rest of the primary solutions, one can assume the solution is $C^{\lfloor A\rfloor}$ and solve for the coefficients up to $\lfloor A\rfloor$ formally. On $x\geq 0$, one can then add a term $b_{A,+} x^A$ where $b_{A,+}$ is a free parameter and recursively solve for a series in $x,x^A$ on $x\geq 0$. On $x\leq 0$, one introduces the term $b_{A,-}(-x)^{A}$ and recursively solves for a series in $x,(-x)^{A}$. Each of these series will necessarily converge, and so for any pair $b_{-,A},b_{+,A}$, one obtains a local primary solution of regularity $C^{\lfloor A\rfloor,A-\lfloor A\rfloor}$.

We can now state the main assumptions on the parameters.
\begin{definition}
  \label{def:admissibility}
  A set of parameters $(d,\alp,\bt,\iota)\in \Z_{\geq 2}\times \R_{>0}\times \R_{>0}\times \{+,-\}$ is \emph{admissible} if the following conditions are met:
  \begin{itemize}
  \item $J_1\geq 0.$ 
  \item $w_0$ defined by \eqref{eq:w_quadratic} with the choice of sign $\iota$, satisfies the physicality condition \eqref{eq:w0-physicality}.
  \item The discriminant $J_2$ (depending on $w_0$) satisfies $J_2> 0$.
  \item $J_2$ satisfies \eqref{eq:J2_ineq} and $A\notin \Z$.
  \end{itemize}
\end{definition}
The preceding discussion can then be summarized in the following Lemma.
\begin{lemma}
\label{lem:background_construction}
  Given admissible data $(d,\alp,\bt,\iota)$, we may associate to it a two-parameter family of primary direction solutions to \eqref{eq:autonomous_ode_sys}, each defined in a neighborhood of $0,$ which satisfy $1-w_0>0$, $\kpp>0,$ and $A>1$. In this family, there is a unique $C^{\infty}$ solution.
\end{lemma}

We state the main assumptions satisfied by the self-similar solutions considered as backgrounds in our results.
\begin{definition}\label{def:admissible_background}
    We say that a self-similar implosion solution is an \textit{admissible regular implosion background} if $(d,\alp,\bt,\iota)$ are admissible parameters in the sense of Definition \ref{def:admissibility} and  $(\hta,\htb)$ is an associated primary direction self-similar solution defined in a neighborhood of $\{x=0\}$.
\end{definition}

\begin{definition}
    \label{def:admissible_background_smooth}
    We say that an admissible regular self-similar implosion solution $(\hta,\htb)$ is an \textit{admissible smooth implosion background} if $(\hta,\htb)$ is $C^{\infty}$.
\end{definition}

We note in the following result that the parameter range considered in \cite{merleImplosionCompressibleFluid2022d} is admissible, so the results of the present paper apply to the examples constructed in \cite{merleImplosionCompressibleFluid2022d}.
\begin{proposition}
  \label{prop:mrrs}
  Define $\bt_+(d,\alp)$ and $\bt_*(d,\alp)$ by
  \begin{equation*}
    \frac{1}{\bt_+}=1+\frac{d-1}{(1+\sqrt{\frac{1}{\alp}})^2},\qquad \frac{1}{\bt_*}=\frac{d+\frac{1}{\alp}}{\frac{1}{\alp}+\sqrt{d}}.
  \end{equation*}
  If $\bt$ satisfies
  \begin{equation*}
    \begin{cases}
      \bt^*<\bt<1&\frac{1}{\alp}<d,\\
      \bt_+<\bt<\bt_*&\frac{1}{\alp}>d,
    \end{cases}
  \end{equation*}
  and $A\notin \Z$, then the data $(d,\alp,\bt,-)$ is admissible. In the language of \cite{merleImplosionCompressibleFluid2022d}, the solution described in the Lemma~\ref{lem:background_construction} is a local solution which passes through a $P2$ point in the $c_-$ direction.
\end{proposition}
\begin{proof}
  See Lemma 2.9 of \cite{merleImplosionCompressibleFluid2022d}.
\end{proof}

\begin{remark}
    The solutions constructed in \cite{merleImplosionCompressibleFluid2022d} have $A\to \infty$ as $\bt_n\to \bt_+$ for $d<\alpha^{-1}$ and $\bt_n\to \bt_*$ for $d>\alpha^{-1}$. For the results in this work, $A$ can be arbitrarily large.
\end{remark}
\subsection{Perturbation Equations}
\label{sec:nonlinear-equations}

The goal of this section is to introduce the main pair of unknowns $(\ap,\bp)$ and derive their evolution equations. We consider \eqref{eq:renormalized_riemann_variables} and substitute
\begin{equation}
  \label{eq:perturbation_def}
  \mra=\hta+\ap,\qquad \mrb=\htb+\bp,
\end{equation}
where $\hta$, $\htb$ is a given self-similar solution. We consider $\htxi,\htc,\httht_1,\httht_2$ to be the associated quantities to this background solution. 

For later use, we note that Section~\ref{sec:append-self-simil} implies that we have the expansions as $x\to 0$
\begin{equation}
  \label{eq:background_taylor_expansions}
  \htxi(x)=1-w_0+O(x),\qquad \htc(x)=\frac{1-w_0}{\alp}+O(x).
\end{equation}
We introduce the notation
\begin{equation*}
  S_1(x)=\htxi-\alp\htc,\qquad S_2(x)=\htxi+\alp\htc.
\end{equation*}
Using \eqref{eq:background_taylor_expansions}, we obtain the following expansions as $x\to 0$
\begin{equation}
  \label{eq:S_taylor_expansions}
  S_1(x)=\kpp x+o(x),\qquad S_2(x)=2(1-w_0)+o(x).
\end{equation}
Using these in \eqref{eq:renormalized_commutator}, it follows that
\begin{equation}
  \label{eq:zeroth_order_commutator}    
  \httht_1(x)=-\kpp+o(1),\qquad \httht_2(x)=o(1).
\end{equation}
We define the following quantities
\begin{equation*}
  \xip=\frac{\ap+\bp}{2},\qquad \cp=\frac{\ap-\bp}{2},
\end{equation*}
which implies that
\begin{equation*}
  \mrxi=\htxi+\xip,\qquad \mrc= \htc+\cp.
\end{equation*}
Using the notation introduced above, we can rewrite the vector fields $\mrT$ and $\mrX$ as follows
\begin{equation*}
  \mrT=\p_{\tau}+(\htxi+\xip)\p_x,\qquad \mrX=\alp(\htc+\cp)\p_x.
\end{equation*}
As a result, the null vector fields are given by
\begin{equation}
  \label{eq:perturbation_vector_fields}
  \mrulL=\p_{\tau}+\left(S_1(x)+\xip-\alp \cp\right)\partial_x,\qquad \mrL=\p_{\tau}+\left(S_2(x)+\xip+\alp \cp\right)\partial_x.
\end{equation}
One obtains the following system for $(\ap,\bp)$, which we refer to as the \emph{perturbation equations}:
\begin{equation}
  \label{eq:perturbation_eqn}
  \begin{cases}
    \mrulL \bp+P_{11}\bp+P_{12}\ap+\calG(\ap,\bp)=0,\\
    \mrL \ap+P_{21}\bp+P_{22}\ap+\calF(\ap,\bp)=0,
  \end{cases}
\end{equation}
where $P_{ij}$ are defined by \eqref{eq:P_entries}, and the quadratic nonlinearities $\calF,\calG$ are defined as
\begin{equation*}
  \calG=\left(-\frac{\alp-1}{2}\ap+\frac{\alp+1}{2}\bp\right)\bp-\frac{(d-1)\alp}{4}\left(\ap^2-\bp^2\right),
\end{equation*}
\begin{equation*}
  \calF=\left(\frac{\alp+1}{2}\ap-\frac{\alp-1}{2}\bp\right)\ap+\frac{(d-1)\alp}{4}\left(\ap^2-\bp^2\right).
\end{equation*}
We note that by \eqref{eq:renormalized_commutator} we also have that
\begin{equation}
  \label{eq:tht_decomposition}
  \mrtht_1=\frac{\mrL(\alp \cp-\xip)-(S_2+\xip+\alp \cp)\p_xS_1}{-S_1+S_2+2\alp \cp},\qquad \mrtht_2=\frac{\mrulL(\alp \cp+\xip)+(S_1+\xip-\alp \cp)\partial_xS_2}{-S_1+S_2+2\alp \cp}.
\end{equation}
We remark that at this point, one can state and prove Theorem \ref{thm:1}.

\subsection{Double Null Formulation}
\label{sec:double-null-form}
In this section we write the system \eqref{eq:perturbation_eqn} using double null coordinates. This coordinate system is convenient for proving Theorem~\ref{thm:2} and \ref{thm:3}.

We begin with the equations \eqref{eq:perturbation_vector_fields}-\eqref{eq:perturbation_eqn}. We define null coordinates $u,v$ by the following transport equations
\begin{equation}
  \label{eq:double_null_def}
  \mrulL v=0,\qquad \mrL u=0.
\end{equation}
We then define inverse foliation densities as
\begin{equation*}
  \mrmu^{-1}:=2\mrT u,\qquad \mrnu^{-1}:=2\mrT v.
\end{equation*}
We rewrite the above vector fields with respect to the double null coordinates $(u,v)$
\begin{equation}
\label{eq:transport_renormalizations}
  \mrulL=\mrmu^{-1}\p_u,\qquad \mrL=\mrnu^{-1}\p_{v},
\end{equation}
\begin{equation*}
  \mrX=\frac{1}{2}(\mrnu^{-1}\p_v-\mrmu^{-1}\p_u),\qquad \mrT=\frac{1}{2}(\mrnu^{-1}\p_v+\mrmu^{-1}\p_u).
\end{equation*}

We commute the transport equations \eqref{eq:double_null_def} with $T$ to obtain transport equations for $\mrmu,\mrnu$. Using the identities in \eqref{eq:renormalized_commutator}, we have
\begin{equation*}
  \p_v\mrmu=\p_{u}\mrnu=-\mrmu\mrnu\mrtht.
\end{equation*}

We further compute that the coordinates $\mrtau,\mrx$ satisfy
\begin{equation}
  \label{eq:coordinate_evolution}
  \p_u\mrtau=\mrmu,\qquad \p_{v}\mrtau=\mrnu,\qquad \p_{u}\mrx=(S_1+\xip-\alp \cp)\mrmu,\qquad \p_{v}\mrx=(S_2+\xip+\alp \cp)\mrnu.
\end{equation}

We recall \eqref{eq:tht_decomposition} and we define renormalizations for $\mrtht_1$ and $\mrtht_2$
\begin{align*}
  \mrTht_1&:=\mrnu\mrtht_1=\frac{\p_v(\alp \cp-\xip)-\mrnu(S_2+\xip+\alp \cp)\p_x S_1}{-S_1+S_2+2\alp \cp},\\
  \mrTht_2&:=\mrmu\mrtht_2=\frac{\p_u(\alp \cp+\xip)+\mrmu(S_1+\xip-\alp \cp)\p_xS_2}{-S_1+S_2+2\alp \cp}.
\end{align*}
We then can rewrite \eqref{eq:perturbation_eqn} as the following system of \emph{null perturbation equations}
\begin{equation}
  \label{eq:double_null}
  \begin{cases}
    \p_{u}\bp+\mrmu \left(P_{11}\bp+P_{12}\ap+\calG(\ap,\bp)\right)=0,\\
    \p_{v}\ap+\mrnu \left(P_{21}\bp+P_{22}\ap+\calF(\ap,\bp)\right)=0,\\
    \p_{u}\mrnu+\mrmu\mrTht_1+\mrnu\mrTht_2=0,\\
    \p_{v}\mrmu+\mrmu\mrTht_1+\mrnu\mrTht_2=0,\\
    \p_{v}\mrx=(S_2+\xip+\alp \cp)\mrnu.
  \end{cases}
\end{equation}
We remark that the functions $P_{ij},S_i,$ and $\p_x S_i$ which appear in the equations \eqref{eq:double_null} are fixed functions of $\mrx$, as defined in \eqref{eq:P_entries} and Section~\ref{sec:nonlinear-equations}. Thus, we include the fifth equation in \eqref{eq:double_null}, in order to obtain a closed system. We also note that the equation for $\p_{u}\mrx$ in \eqref{eq:coordinate_evolution} holds initially and is propagated by the equations \eqref{eq:double_null}.

Finally, we also mention that when $\ap=\bp=0$, the perturbation is trivial, so $\mra$ and $\mrb$ simply recover their background values $\hta$ and $\htb.$ We make the notation convention that $\httau(u,v),$ $\htx(u,v),$ $\htTht_i(u,v),$ $\htP_{ij}(u,v),$ $\htmu(u,v),$ and $\htnu(u,v)$ correspond to this solution. These functions and their relevant asymptotics as $v\to 0$ are constructed in the following section.
\subsection{Background Null Coordinates}
\label{sec:append-null-coord}
We consider an admissible regular self-similar solution, as in Definition \ref{def:admissible_background}. We will consider the case when $(\hta,\htb)\in C^{2}$. This includes the case $A>2$, as in Theorem \ref{thm:2}. This also covers the case of $A>1$ and the $C^{\infty}$ primary direction solution, as in Theorem \ref{thm:3}, in which case all of the regularities in this section may be replaced by $C^{\infty}$. In this section, we derive asymptotics for the corresponding double null coordinates which satisfy a certain gauge condition. We will normalize the gauge so that $(\tau,x)=(0,0)$ corresponds to $(u,v)=(0,0)$. We then impose
\begin{equation}
\label{eq:background_gauges}
    v|_{u=0}=\tau,\qquad u|_{v=0}=\tau.
\end{equation}
We comment that by the identity \eqref{eq:transport_renormalizations}, this is equivalent to
\begin{equation*}
    \htnu|_{u=0}=\htmu|_{v=0}=1.
\end{equation*}
We comment that although the background solution is independent of $\tau,$ the coordinates $(u,v)$ depend on both $\tau$ and $x$. In the remainder of this section, we derive the asymptotics of these null coordinates near $x=0$. We begin by introducing the following spatial coordinate
\begin{equation*}
  x^*(x)=\int_{0}^{x}(\htxi+\alp \htc)^{-1}dx'.
\end{equation*}
As a result, it follows that
\begin{equation}
  \label{eq:u_background_def}
  u=\tau-x^*
\end{equation}
is a null coordinate which satisfies the gauge condition \eqref{eq:background_gauges}. 

We define another spatial coordinate $\tlx$, by solving the following singular ODE:
\begin{equation*}
  (\htxi-\alp\htc) \frac{\mathrm{d}\tlx}{\mathrm{d} x}-\kpp\tlx=0,\qquad \tlx(0)=0,\qquad \tlx'(0)=1.
\end{equation*}
We note that, since this is a singular ODE, one must prescribe $\tlx'(0)$ and then there is a unique $C^2$ solution locally.

It follows that $\tlv:=e^{-\kpp\tau}\tlx$ is a null coordinate. We remark that by \eqref{eq:u_background_def}, as well as \eqref{eq:background_taylor_expansions}, we have that
\begin{equation*}
  \tlv|_{u=0}=2(1-w_0)\tau+O(\tau^2).
\end{equation*}
In particular, there exists a $\Psi$ defined locally with
\begin{equation*}
  \Psi(0)=0,\qquad \Psi'(0)=\frac{1}{2(1-w_0)},
\end{equation*}
with the property that the rescaled coordinate
\begin{equation}
  \label{eq:v_background_def}
  v:=\Psi(\tlv)
\end{equation}
is a null coordinate satisfying the gauge condition \eqref{eq:background_gauges}.

For the rest of the section, we express the coordinates $(\tau,x)$ in terms of the null coordinates $(u,v).$ We then derive asymptotic expansions for the background quantities as $v\rightarrow0$.

Firstly, we note that for
\begin{equation*}
  z_1:=x,\qquad z_2:=xe^{-\kpp \tau},
\end{equation*}
we have that locally
\begin{equation*}
  u=\frac{1}{\kpp}\ln(z_1/z_2)-\frac{z_1}{2(1-w_0)}F_1(z_1,z_2),\qquad F_1\in C^{1},\qquad F_1(0)=1,
\end{equation*}
\begin{equation*}
  v=\frac{z_2}{2(1-w_0)}F_2(z_1,z_2),\qquad F_2\in C^{1},\qquad F_2(0)=1.
\end{equation*}
In the other direction, we define
\begin{equation}
  \label{eq:w_def}
  w_1=2(1-w_0)e^{\kpp u}v,\qquad w_2=2(1-w_0)v.
\end{equation}
We now consider the $\tau,x$ coordinates as functions of $u,v$. We denote these functions as $\httau,\htx$, to emphasize that these expressions hold for the background solution. We have locally
\begin{equation}
  \label{eq:x_asymptotics}
  \htx=w_1G_1(w_1,w_2),\qquad G_1\in C^{1},\qquad G_1(0)=1,
\end{equation}
\begin{equation*}
  \httau=\frac{1}{\kpp}\ln(w_1/w_2)+\frac{w_1}{2(1-w_0)}G_2(w_1,w_2),\qquad G_2\in C^{1},\qquad G_2(0)=1.
\end{equation*}
We comment that as corollaries of the above expansions, we have locally
\begin{equation*}
  z_2=w_2\cdot G_3(w_1,w_2),\qquad G_3\in C^{1},\qquad G_3(0)=1,
\end{equation*}
\begin{equation*}
  w_1=z_1\cdot F_3(z_1,z_2),\qquad F_3\in C^{1},\qquad F_3(0)=1.
\end{equation*}
In particular $(z_1,z_2)\to (w_1,w_2)$ is a regular change of variables at $0,$ whose linearization is the identity.

We comment that in particular, we have the expansions
\begin{equation*}
  \httau=u+e^{\kpp u}v+O((e^{\kpp u}v)^2),\qquad \htx=2(1-w_0)e^{\kpp u}v+O((e^{\kpp u}v)^2).
\end{equation*}
We also comment that the derivatives of the null coordinates are given by
\begin{equation*}
  \begin{pmatrix}
    \pdv[\httau]{u}&\pdv[\httau]{v}\\
    \pdv[\htx]{u}&\pdv[\htx]{v}
  \end{pmatrix}=
  \begin{pmatrix}
    \htmu & \htnu\\
    -\htmu(\alp\htc-\htxi) & \htnu(\alp\htc+\htxi)
  \end{pmatrix}.
\end{equation*}
As a result, it follows that we have the expansions
\begin{equation}
  \label{eq:munuasymptots}
  \htmu=1+O(e^{\kpp u}v)
  ,\qquad \htnu=e^{\kpp u}+O(e^{2\kpp u}v).
\end{equation}
Finally, one can then express $\htc,\htxi$ in terms of $u,v$. We have by \eqref{eq:S_taylor_expansions}
\begin{equation*}
  \htxi-\alp\htc=2(1-w_0)\kpp e^{\kpp u}v+O(e^{2\kpp u}v^2).
\end{equation*}
In particular, we also obtain the expansions
\begin{equation}
  \label{eq:Tht_asymptotics}
  \htTht_1=-\kpp e^{\kpp u}+O(e^{2\kpp u}v),\qquad \htTht_2=O(e^{\kpp u}v).
\end{equation}
\subsection{Gauge for Double Null Perturbation System}\label{sec:pertubation_gauge}
We now return to the null perturbation equations \eqref{eq:double_null}. In our second main result, we will consider perturbations which vanish for $x\leq 0$. On the other hand, for our third main result, we will consider perturbations which do not necessarily vanish for $x\leq 0$. In this section, we fix the choice of double null coordinates for both such perturbations. 

We consider an outgoing $\mrL$ characteristic curve which emanates from $p_-:=(\tau_-,x_-)$ and passes through $(\tau,x)=(0,0)$. In the following discussion, we consider $(\tau_-,x_-)$ to be parameters, where $x_-\leq 0$ is sufficiently close to $0$. We define $\dlt'\geq 0$ implicitly by the relation
\begin{equation}
    \label{eq:x_dlt_conversion}
    \htx(0,-\dlt')=x_-.
\end{equation}
We normalize the null coordinates by letting the point $(\tau,x)=(\tau_-,x_-)$ correspond to $(u,v)=(0,-\dlt')$.

We fix the double null gauge by imposing the following conditions, as long as the solutions to the resulting ODEs are defined
\begin{equation}\label{eq:densities_v=-dlt}
    \mrnu|_{u=0}=1,\qquad \mrmu|_{v=-\dlt'}=\htmu|_{v=-\dlt'}.
\end{equation} 

The first ODE can be solved explicitly to give
\begin{equation*}
    v|_{u=0}=-\dlt'+(\tau-\tau_-).
\end{equation*}

We note that by \eqref{eq:coordinate_evolution}, if $\ap=\bp=0$ along $\{v=-\dlt'\}$, then by \eqref{eq:densities_v=-dlt} we have that
\begin{equation*}
    \mrx|_{v=-\dlt'}=\htx|_{v=-\dlt'},\qquad \mrtau|_{v=-\dlt'}=\httau|_{v=-\dlt'}+\tau_-+\delta'.
\end{equation*}

We conclude this section by briefly explaining the relation between the choice of $p_-$ and the perturbations considered. For our second main result, we will consider exterior perturbations, so we set $p_-=(0,0)$, which corresponds to $\delta'=0.$ On the other hand, for our third main result, we allow for interior perturbations as well, so we will take $\delta'>0.$

\section{Main Theorems}
\label{sec:shock-form-results}
We present three results on shock formation for the perturbation system \eqref{eq:perturbation_eqn}. This is a $(1+1)$-dimensional symmetric hyperbolic system and we will always work with strong solutions of regularity $C^{k,1}$ for some $k\geq 0$. The system has finite speed of propagation and local well-posedness in $C^{k,1}$ is classical. We comment that by ``the solution'' we will always refer to a maximal globally hyperbolic development (MGHD) of data. In particular, we do not discuss discontinuous weak solutions.

\subsection{Exterior Perturbations of Regularity $C^{N,1}$}
In our first two results, we will consider perturbations with finite regularity.

We recall that the \emph{interior} and \emph{exterior} regions of spacetime are given by $x\leq 0$ and $x\geq 0$ respectively. For the results of this section, we will assume the trivial interior solution: $\ap=\bp=0$ for $x\leq 0$ which leaves the background unchanged in the interior region.

In order to state and prove these theorems is useful to work with \emph{exterior solutions}. These are defined to be solutions to \eqref{eq:perturbation_eqn} on $x\geq 0$. It is assumed that $(\ap,\bp)$ both vanish identically on left boundary $x=0$. At this boundary, $C^{N}$ regularity is to be understood in the sense of one-sided derivatives. Given an exterior solution, it may be glued along $x=0$ to the interior solution. This \emph{glued solution} will typically have a jump in the derivative and hence be of regularity $C^{0,1}$. Using the equations, one can see that the obstructions to higher regularity are given by $(\mrL)^{k}\bp|_{x=0}$. This is summarized in the following lemma.
\begin{lemma}
    \label{lem:gluing}
    Let $(\ap,\bp)$ be an exterior solution of regularity $C^{N}$. Suppose that
    \begin{equation*}
        \mrL \bp|_{x=0}=\cdots=(\mrL)^{(N-1)} \bp|_{x=0}=0,\qquad (\mrL)^{(N)}\bp|_{x=O}\neq 0.
    \end{equation*}
    Then, $(\ap,\bp)$ can be glued to the trivial interior solution to obtain a glued solution which is piecewise $C^{N}$ with a jump in the $N$th derivative at $x=0$, and hence in $C^{N-1,1}\setminus C^N$.
\end{lemma}
\begin{proof}
    Under the above assumptions, by commuting the $\ap$ equation of \eqref{eq:perturbation_eqn} with $\mrL$ a total of $N-1$ times, one can inductively show $\mrL^{(k)}\ap|_{x=0}=0$ for $1\leq k\leq N$.
\end{proof}

\begin{remark}
    In the statement and proofs of the first two theorems, $C^{N,1}$ refers to the regularity of the exterior perturbation. See for instance, the title of this subsection, the statements of Definition~\ref{def:char_data}, Theorem~\ref{thm:1} and Theorem~\ref{thm:2}. On the other hand, for summarizing these results, it is helpful to think of the glued solution which is of regularity $C^{N-1,1}$, see for instance the theorem statements in the introduction and the statement of Corollary \ref{cor:instability_corollary}. 
\end{remark}

It will also be convenient to consider the exterior solution as the development of characteristic initial data. The characteristics will emanate towards the future from
\begin{equation*}
  O:=(\tau=0,x=0).
\end{equation*}
In this case, data for $\ap$ will be imposed along the ingoing characteristic, which we denote by
\begin{equation}
  \label{eq:ulcalC_def}
  \ulcalC_O:=\left\{(\tau,x)\colon 0\leq \tau<\tau_0,x=0\right\}.
\end{equation}
Here $0<\tau_0\leq \infty$ is a parameter, and we will impose $\ap|_{\ulcalC_O}=0$. We will also always impose $\bp|_{O}=0$. These assumptions ensure that $\ulcalC_O$ is an ingoing $\mrulL$-characteristic and $\bp|_{\ulcalC_O}=0$. Data for $\bp$ will be imposed on an outgoing $\mrL$-characteristic, which we will denote by $\calC_O$. The precise assumptions on the data are given in the following definition.
\begin{definition}
\label{def:char_data}
    For the perturbation equations \eqref{eq:perturbation_eqn} and $N\geq 1$, we define \emph{$C^{N,1}$ exterior characteristic perturbation data} to be a $C^{N,1}$ curve
    $\calC_O$ emanating from $O$ and a function $\brb\in C^{N,1}(\calC_O)$ such that $\brb|_{O}=0$. We assume $\ulcalC_O$ is given by \eqref{eq:ulcalC_def} with $\tau_0=\infty$ and $\ap|_{\ulcalC_O}=0$ (in particular $\ulcalC_O$ is an ingoing characteristic). We impose that $\mrL$ is tangent to $\calC_O$ (i.e. $\calC_O$ is outgoing characteristic).
\end{definition} 

We comment that the assumption of characteristic data is not essential. To construct data for these theorems, one could start with small data for $(\ap,\bp)$ at $\tau=0$ and then (locally) evolve up to $\calC_O$.

We also comment that in the case of a regular, but not $C^{\infty}$, background and $N=\lfloor A\rfloor$, the regularity assumed above for the perturbation is higher than the regularity of the background, $C^{\lfloor A\rfloor,1}$ as opposed to $C^{\lfloor A\rfloor,A-\lfloor A\rfloor}$. This does not impact the proof and we comment that one could replace $C^{N,1}$ with $C^{N,\alp'}$ for any $0<\alp'<1$ in Definition \ref{def:char_data} and obtain similar results.

We also remark that the results are stated as blow-up criteria. Nonetheless, one can take the initial data perturbation to be arbitrarily small in the relevant $C^{k,1}$ norm. Thus, one should think of these results as describing instability of the background solutions, which is stated as Corollary~\ref{cor:instability_corollary}.

We comment that the results are local in nature. By finite speed of propagation, these results imply that any strong solutions obeying the hypotheses must blow up in the classical sense. If one considers weak solutions, an incoming shock could in principle cover the singularities we construct in this section. However, we note that the solutions constructed later in Theorem~\ref{thm:3} extend up to $\{t=0\}$ in the outer characteristic region, ruling out incoming shocks from this region.

We now state the first result, which is the detailed formulation of Theorem~\ref{thm:1_intro}.
\begin{theorem}
  \label{thm:1}
  Let $(\hta,\htb)$ be an admissible regular implosion background as in Definition~\ref{def:admissible_background}. We consider the perturbation equations \eqref{eq:perturbation_eqn} with $C^{1,1}$ exterior characteristic perturbation data as in Definition~\ref{def:char_data}. Suppose that, as a one-sided derivative along $\calC_O$,
  \begin{equation*}
    \mrL \brb|_{O}<0.
  \end{equation*}
  Define a constant $\lmb$ by
  \begin{equation}
  \label{eq:lmb_def}
  \lmb:=\frac{\alp+1}{4(1-w_0)(A-1)\kpp}
\end{equation}
  and a time $\tau_S$ implicitly by
  \begin{equation}
    \label{eq:blowup_timescale}
    e^{\kpp (A-1)\tau_S}=1-\frac{1}{\lmb\cdot \mrL \brb|_{O}}.
  \end{equation}
  It follows that $\bp|_{x=0}=0$ and the one-sided derivative $\mrL \bp|_{x=0}$ blows up at $\tau=\tau_S$.
\end{theorem}
\begin{proof}
  We obtain from the equations \eqref{eq:perturbation_eqn} that $\bp|_{x=0}=\mrulL \bp|_{x=0}=0,$ $\mrL \ap|_{x=0}=0$, and $\xip|_{x=0}=\cp|_{x=0}=0$. In particular, we also get that $\mrulL|_{x=0}=\partial_{\tau}.$ We note that as a consequence, we can compute using \eqref{eq:tht_decomposition}
  \begin{equation}
    \label{eq:tht_ode_compute}
    \mrtht_1|_{x=0}=-\kpp-\frac{\alp+1}{4(1-w_0)}\mrL \bp|_{x=0},\qquad \mrtht_2|_{x=0}=0.
  \end{equation}
  We now commute the first equation in \eqref{eq:perturbation_eqn} with $\mrL$. We then obtain the equation
  \begin{equation*}
    \mrulL(\mrL \bp)+[\mrL,\mrulL]\bp+[\mrL,P_{11}]\bp+P_{11}\mrL \bp+[\mrL,P_{12}]\ap+P_{12}\mrL \ap+\mrL(\calG(\ap,\bp))=0.
  \end{equation*}
  Using \eqref{eq:renormalized_commutator} and \eqref{eq:tht_ode_compute}, we have
  \begin{align*}
    [\mrT,\mrX]\bp|_{x=0}&=
    -\frac{\kpp}{2}\mrL \bp|_{x=0}-\frac{\alp+1}{8(1-w_0)}\left(\mrL \bp|_{x=0}\right)^2.
\end{align*}

We denote $b_1=\mrL \bp|_{x=0}$. Using the above equations, we obtain the following ODE along $\{x=0\}$
\begin{equation}
  \label{eq:b_ricatti}
  \frac{db_1}{d\tau}-\kpp(A-1)b_1+\frac{\alp+1}{4(1-w_0)}b_1^2=0.
\end{equation}
As a comment, the quadratic term has a positive coefficient. As a consequence, one sees a finite-time blow-up if $b_1(0)<0$. Solving the ODE yields the timescale \eqref{eq:blowup_timescale}.
\end{proof}
\begin{remark}
\label{rem:op}
  We remark that similar arguments were previously considered by Ori and Piran in their investigations on the stability of self-similar gravitational collapse profiles for the Euler-Poisson system, \cite{oriSimpleStabilityCriterion1988}. 
  Ori-Piran's stability criterion, using the language of this paper, is $\kpp>0$ and $A<1$. They show that in the isothermal case, $\gmm=1$, this stability criterion is satisfied for the Larson-Penston and Hunter self-similar solutions, and conclude stability against $C^{0,1}$ exterior perturbations.
\end{remark}

We now state the second result.
\begin{theorem}
  \label{thm:2} 
  Consider an admissible regular implosion background as in Definition~\ref{def:admissible_background}. Fix $2\leq N<A,$ and consider the perturbation equations \eqref{eq:perturbation_eqn} with $C^{N,1}$ exterior characteristic perturbation data as defined in Definition \ref{def:char_data}. Suppose that the perturbation data satisfies
  \begin{equation*}
    \mrL\brb|_{O}=\cdots=(\mrL)^{(N-1)}\brb|_{O}=0,\qquad (\mrL)^{(N)}\brb|_{O}<0.
  \end{equation*}
  In this case, in the domain of dependence of sufficiently small subcurves $\calC_O'\subset \calC_O$ emanating from $O$, the solution of \eqref{eq:perturbation_eqn} blows up in $C^{0,1}$ while remaining bounded in $L^{\infty}$.
\end{theorem}

As a corollary of Theorem \ref{thm:1} and \ref{thm:2}, we have the following instability result, where we recall the notion of glued solution introduced in the discussion before Lemma~\ref{lem:gluing}.
\begin{corollary}
\label{cor:instability_corollary}
    Consider an admissible regular implosion background as in Definition~\ref{def:admissible_background}. Then for $0\leq k\leq \lfloor A\rfloor -1$, the zero perturbation has $\lfloor A\rfloor-k$ unstable directions in the space of strong, glued solutions, equipped with the $C^{k,1}$ topology.
\end{corollary}
\begin{proof}
    This follows from Theorem \ref{thm:1} and \ref{thm:2}, where we take the perturbations to be small and we observe that the perturbations for different $N$ are linearly independent.
\end{proof}
To prove Theorem \ref{thm:2}, we will work in the double null formulation given by the null perturbation equations \eqref{eq:double_null}, with the gauge conditions in Subsection~\ref{sec:pertubation_gauge} for $\delta'=0$ (as we only consider exterior perturbations in this case). We recall that in these coordinates, $O$ corresponds to $(u,v)=(0,0)$. We comment that when working in double null coordinates, we will always have for some $\Delta>0$
\begin{equation}
  \label{eq:calC_def_null}
  \ulcalC_O:=\left\{(u,v)\colon v=0,\ 0\leq u<\infty\right\},\qquad \calC_O:=\left\{(u,v)\colon u=0,\ 0\leq v\leq \Delta\right\}.
\end{equation}

We will deduce Theorem \ref{thm:2} at the end of this section as a consequence of the following result for the system of null perturbation equations.
\begin{lemma}
  \label{lem:null_shock_result}
  Consider an admissible regular implosion background as in Definition~\ref{def:admissible_background}. We define $\ulcalC_O$, $\calC_O$ by \eqref{eq:calC_def_null}, and we consider characteristic data $\ap|_{v=0}=0$ and $\bp|_{u=0}=\brb(v)$ for the null perturbation system \eqref{eq:double_null}. We fix an integer $2\leq N<A$, and we assume that $\brb\in C^{N,1}(\calC_O)$ satisfies
  \begin{align*}
    \brb(0)=\p_v \brb(0)=\cdots=\p_v^{(N-1)}\brb(0)=0,\qquad \p_v^{(N)}\brb(0)<0.
  \end{align*}
  We impose the gauge conditions
  \begin{equation*}
    \mrmu|_{\ulcalC_O}=1,\qquad \mrnu|_{\calC_O}=1,\qquad \mrx|_{\ulcalC_O}=0.
  \end{equation*}
  In this case, for $\dlt>0$ sufficiently small, there exist explicit $(u_-,u_+)$ depending on the data so that if we define
  \begin{equation}
    \label{eq:calA_def}
    \calA:=[0,u_+)\times [0,\dlt],
  \end{equation}
  then the solution of \eqref{eq:double_null} exists in $\calA$, satisfies uniform $C^1$ bounds and satisfies $\mrnu(u_S,\dlt)=0$ for some $u_S\in (u_-,u_+)$.
  Moreover, throughout the domain $\calA$, we have $\mrmu>0,\ \p_v\bp<0$.
\end{lemma}
\begin{remark}
  \label{rem:mghd}
  While we started with the perturbation equations \eqref{eq:perturbation_eqn} and derived the null perturbation system \eqref{eq:double_null}, solutions to the latter system remain regular even if $\mrmu$ or $\mrnu$ vanish. As a consequence, given corresponding data for the two systems, an MGHD for the null perturbation equations \eqref{eq:double_null} can contain the boundary of an MGHD for the perturbation equations \eqref{eq:perturbation_eqn}.
  Indeed, Lemma \ref{lem:null_shock_result} constructs an MGHD $\calA$ of \eqref{eq:double_null}, for which $\mrnu$ vanishes. To obtain an MGHD of \eqref{eq:perturbation_eqn}, we consider
  \begin{equation*}
    \calA_+:=\left\{(u,v)\in \calA\colon \mrmu,\mrnu>0\right\}.
  \end{equation*}
  We then let $\tilde{\calA}$ be the connected component of $\calA_+,$ which contains the initial data characteristics. We denote the upper boundary of this set by
  \begin{equation*}
    \calS:=\partial_{+} \tilde{\calA}\subset \calA.
  \end{equation*}
  Finally, we consider
  \begin{equation*}
    \calA_0:=\calA\setminus J^+(\calS),
  \end{equation*}
  where $J^+(\calS)$ is the causal future defined as
  \begin{equation*}
    J^+(\calS):=\left\{(u,v)\in \calA\colon \exists\ u'\leq u,v'\leq v,(u',v')\in \calS\right\}.
  \end{equation*}
  It follows that $(u,v)\to(\tau,x)$ is a diffeomorphism on $\calA_0$, so $\calA_0$ is an MGHD for \eqref{eq:perturbation_eqn}. For later reference, we denote the upper boundary of $\calA_0$ by
  \begin{equation*}
    \calS_0:=\partial_+ \calA_0.
  \end{equation*}
  We refer the reader to the recent work \cite{abbrescia2025quasilinearwavesupersonicshock} for a discussion regarding MGHDs and the structure of their boundaries.
\end{remark}
We now prove Lemma \ref{lem:null_shock_result}.
\begin{proof}
  The result will be proved using a bootstrap argument. We will work on the domain $\calA$ as defined by \eqref{eq:calA_def}, with suitable $u_+$ and $\dlt$ to be determined in the course of the proof. The domain is illustrated in Figure \ref{fig:characteristic_domain} in $(\tau,x)$ and $(u,v)$ coordinates. The bold black curve corresponds to $\calS$ and the blue region corresponds to $\calA_0$. The red arrows indicate the causality of the acoustical metric. 

  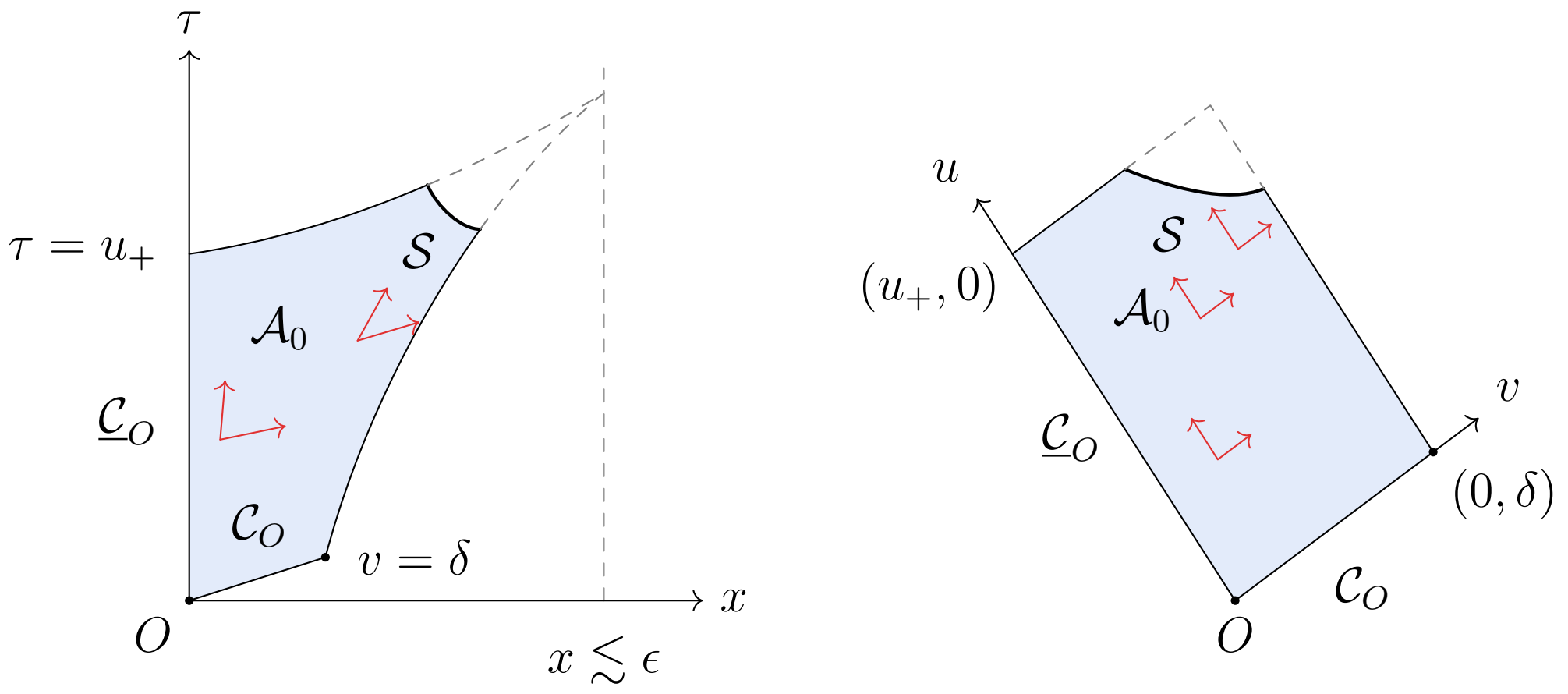
\begin{figure}[h]
    \centering
\begin{tikzpicture}[scale=1.05,line width=0.4pt,
    line cap=round,line join=round]
    \definecolor{domainblue}{RGB}{40,103,224}
    \definecolor{domainred}{RGB}{225,48,48}
    \tikzset{
        singular/.style={black,line width=0.8pt},
        continuation/.style={dashed,black!50},
        causal/.style={->,domainred}
    }

    \begin{scope}
        \coordinate (LO) at (0,0);
        \coordinate (LI) at (1.10,0.35);
        \coordinate (LT) at (0,2.80);
        \coordinate (LF) at (3.35,4.10);
        \coordinate (LP) at (1.923756,3.357074);
        \coordinate (LR) at (2.350611,2.998832);

        \fill[domainblue,fill opacity=0.13]
            (LO) -- (LI)
            .. controls (1.343200,1.278000)
                    and (1.799392,2.226480) .. (LR)
            .. controls (2.200611,2.998832)
                    and (2.003756,3.177074) .. (LP)
            .. controls (1.276000,3.079792)
                    and (0.638000,2.892800) .. (LT)
            -- cycle;

        \draw (LO) -- (LI);
        \draw (LI) .. controls (1.343200,1.278000)
                           and (1.799392,2.226480) .. (LR);
        \draw (LT) .. controls (0.638000,2.892800)
                           and (1.276000,3.079792) .. (LP);
        \draw[singular]
            (LP) .. controls (2.003756,3.177074)
                        and (2.200611,2.998832) .. (LR);

        \draw[continuation]
            (LP) .. controls (2.392820,3.557864)
                        and (2.867000,3.806000) .. (LF);
        \draw[continuation]
            (LR) .. controls (2.660672,3.433280)
                        and (3.000800,3.812000) .. (LF);

        \draw[dashed,black!40] (3.35,0) -- (3.35,4.30);
        \node[anchor=north] at (3.35,-0.15) {$x\lesssim\epsilon$};

        \draw[causal] (0.25,1.30) -- (0.29,1.78);
        \draw[causal] (0.25,1.30) -- (0.78,1.41);
        \draw[causal] (1.36,2.10) -- (1.60,2.53);
        \draw[causal] (1.36,2.10) -- (1.86,2.25);

        \draw[->] (LO) -- (4.15,0) node[right] {$x$};
        \draw[->] (LO) -- (0,4.45) node[above] {$\tau$};
        \fill (LO) circle (1pt);
        \fill (LI) circle (1pt);
        \node[below left] at (LO) {$O$};
        \node[anchor=east] at (-0.12,1.45)
            {$\underline{\mathcal C}_O$};
        \node[anchor=south] at (0.57,0.3) {$\mathcal C_O$};
        \node[anchor=west] at (1.22,0.35) {$v=\delta$};
        \node[anchor=east] at (-0.12,2.80) {$\tau=u_+$};
        \node at (0.73,2.20) {$\mathcal A_0$};
        \node[anchor=north east] at (2.13,3.10) {$\mathcal S$};
    \end{scope}

    \begin{scope}[shift={(8.45,0)},
        x={(-1.80cm,2.80cm)},y={(1.60cm,1.20cm)}]

        \fill[domainblue,fill opacity=0.13]
            (0,0) -- (0,1) -- (0.76,1)
            .. controls (0.78,0.85) and (0.90,0.68) .. (1,0.57)
            -- (1,0) -- cycle;

        \draw (1,0) -- (1,0.57);
        \draw (0,1) -- (0.76,1);
        \draw[singular]
            (1,0.57) .. controls (0.90,0.68)
                            and (0.78,0.85) .. (0.76,1);

        \draw[continuation] (1,0.57) -- (1,1) -- (0.76,1);

        \foreach \u/\v in {0.30/0.25,0.60/0.50,0.68/0.78} {
            \draw[causal] (\u,\v) -- ({\u+0.12},\v);
            \draw[causal] (\u,\v) -- (\u,{\v+0.17});
        }

        \draw[->] (0,0) -- (1.16,0) node[above left] {$u$};
        \draw[->] (0,0) -- (0,1.23) node[above right] {$v$};
        \fill (0,0) circle (1pt);
        \fill (0,1) circle (1pt);
        \node[below] at (0,0) {$O$};
        \node[anchor=east,xshift=-3pt] at (0.47,0)
            {$\underline{\mathcal C}_O$};
        \node[anchor=north west,yshift=-3pt] at (0,0.42)
            {$\mathcal C_O$};
        \node[anchor=east,xshift=1pt,yshift=-8pt] at (1,0) {$(u_+,0)$};
        \node[below right] at (0,1) {$(0,\delta)$};
        \node at (0.70,0.32) {$\mathcal A_0$};
        \node[anchor=north,yshift=-3pt] at (0.90,0.68)
            {$\mathcal S$};
    \end{scope}
\end{tikzpicture}
    \caption{Sketch of the characteristic domain $\mathcal{A}$ for Theorem~\ref{thm:2}.}
    \label{fig:characteristic_domain}
\end{figure}
  
  We consider the $C^{N,1}$ exterior solution to \eqref{eq:double_null} determined by the characteristic data described above. The first two equations in \eqref{eq:double_null} imply that
  \begin{equation*}
    \bp|_{v=0}=\p_v\bp|_{v=0}=\ap|_{v=0}=\p_v\ap|_{v=0}
    =0.
  \end{equation*}
  We define
  \begin{equation*}
    b_N:=\p_v^{(N)}\brb(0).
  \end{equation*}
  By assumption, there exists a Lipschitz constant $L_{N+1}>0$, so that for any $0\leq v\leq\delta$
  \begin{equation}
    \label{eq:b_initialization}
    \left|\p_v \brb(v)-b_N \frac{v^{N-1}}{(N-1)!}\right|\leq L_{N+1} \frac{v^{N}}{N!}.
  \end{equation}
  We define a timescale $u_*$ by
  \begin{equation}
    \label{eq:timescale_def}
    \lmb b_N\frac{\dlt^{N-1}}{(N-1)!}\left(e^{\kpp (A-1)u_*}-1\right)=-1.
  \end{equation}
  The claim we will establish is that, for $\dlt$ sufficiently small, along the characteristic $\{v=\dlt\}$, $\mrnu$ will vanish at $u\approx u_*$.

  We will use ``$\lesssim$'' to denote inequalities with implicit constants which only depend on the background ($\alp,\bt,d,\iota$, and by extension $A,\kpp$). We also consider $\ep>0$ to be a sufficiently small constant depending on the background.

  We will solve the equations \eqref{eq:double_null} for $u\in[0,u_+],$ where $u_+>u_{*}$ and $u_{\pm}$ are defined by
  \begin{equation}
    \label{eq:upper_time}
    \lmb b_N\frac{\dlt^{N-1}}{(N-1)!}(e^{\kpp (A-1)u_{\pm}}-1)=-1\mp\ep^{1-\frac{3}{100}}.
  \end{equation}
  By choosing $\dlt$ smaller if necessary, we may impose the following conditions
  \begin{equation}
    \label{eq:dlt_smallness}
    \dlt\leq \ep^{\frac{A-1}{A-N}}\abs{b_N}^{\frac{1}{A-N}},
  \end{equation}
  \begin{equation}
    \label{eq:dlt_smallness_B}
    \dlt \abs{L_{N+1}}\leq \abs{b_N},
  \end{equation}
  \begin{equation}
    \label{eq:zeroth_term_boundedness}
    \lmb\frac{\dlt^{N-1}}{(N-1)!}\abs{b_N}\leq 1.
  \end{equation}
  Combining \eqref{eq:upper_time} and \eqref{eq:zeroth_term_boundedness} implies that
  \begin{equation}
    \label{eq:timescale_bounded}
    \lmb \abs{b_N}\frac{\dlt^{N-1}}{(N-1)!}e^{\kpp(A-1)u_+}\leq 3.
  \end{equation}
  As a consequence of \eqref{eq:timescale_bounded}, we have
  \begin{equation*}
      \dlt e^{\kpp u_+}\lesssim \dlt^{\frac{A-N}{A-1}}\frac{1}{\abs{b_N}^{\frac{1}{A-1}}}.
  \end{equation*}
  By \eqref{eq:dlt_smallness} and the assumption $N<A$, we have that for all $u\in[0,u_+]$
  \begin{equation}
    \label{eq:x_smallness}
    \dlt e^{\kpp u}\lesssim \ep
  \end{equation}
  Moreover, we recall that according to Section~\ref{sec:append-null-coord}, the background functions can be Taylor expanded in $w_1,w_2$, defined by \eqref{eq:w_def}, which are of order $ve^{\kpp u}$, $v$. By choosing $\ep$ smaller if necessary, we may approximate background functions by their Taylor expansions. 

  We also note that by \eqref{eq:x_asymptotics}, which will hold to leading order, we will have
  \begin{equation*}
    x\approx 2(1-w_0)ve^{\kpp u},
  \end{equation*}
  so $\ep$ measures the maximum $x$-width of the characteristic region (see Figure \ref{fig:characteristic_domain}).
  
  As a consequence of \eqref{eq:b_initialization} and \eqref{eq:dlt_smallness_B}, we have the bounds
  \begin{equation*}
    \abs{\brb}\lesssim \abs{b_N}v^{N},\qquad \abs{\p_v \brb}\lesssim \abs{b_N}v^{N-1}
    .
  \end{equation*}
  We consider the following five bootstrap assumptions:
  \begin{align*}
    \abs{\p_v \bp-e^{\kpp A u}\p_v\brb}&\leq \ep^{-\frac{2}{100}}\abs{b_N} e^{\kpp (A+1) u}v^{N},\\
    \abs{\p_u \ap}&\leq \ep^{-\frac{1}{100}}\abs{b_N} e^{\kpp (A+1) u}v^{N+1},\\
    |\mrnu-\htnu-\lmb (e^{\kpp A u}-e^{\kpp u})\p_v\brb|&\leq \ep^{-\frac{3}{100}}\abs{b_N} e^{\kpp (A+1) u}v^{N},\\
    \abs{\mrmu-\htmu}&\leq \ep^{-\frac{1}{100}}\abs{b_N}e^{\kpp A u}v^N,\\
    \abs{\mrx-\htx}&\leq \ep^{-\frac{1}{100}}\abs{b_N}e^{\kpp A u}v^N.
  \end{align*}
  It is straightforward to check that the bootstrap assumptions hold initially with significantly better constants. As a comment, integrating the above inequalities yields
  \begin{equation*}
    \abs{\bp-e^{\kpp A u}\brb}\lesssim \ep^{-\frac{2}{100}}\abs{b_N} e^{\kpp (A+1) u}v^{N+1},\qquad \abs{\ap}\lesssim \ep^{-\frac{1}{100}}\abs{b_N}e^{\kpp(A+1)u}v^{N+1}.
  \end{equation*}
  We then obtain from the triangle inequality and \eqref{eq:x_smallness} that
  \begin{equation*}
    \abs{\p_v \bp}\lesssim \abs{b_N}e^{\kpp A u}v^{N-1},\qquad
    \abs{\bp}\lesssim \abs{b_N}e^{\kpp A u}v^{N},\qquad
    \abs{\mrnu-\htnu}\lesssim \abs{b_N}e^{\kpp A u}v^{N-1}.
  \end{equation*}
  One can use \eqref{eq:timescale_bounded} and the previous bounds to obtain the following estimate for nonlinear terms, which will be used below
  \begin{equation*}
    \abs{\bp^2}\lesssim e^{\kpp(A+1)u}\abs{b_N}v^{N+1},\qquad \abs{\bp\p_v \bp}\lesssim e^{\kpp(A+1)u}\abs{b_N}v^{N},\qquad
    \left|(\mrnu-\htnu)\bp\right|\lesssim \abs{b_N}e^{\kpp(A+1)u}v^N.
  \end{equation*}
  \emph{Closing the bootstrap argument.}
    Using the bootstrap assumption on $\mrx$ in \eqref{eq:A_property_a} and \eqref{eq:S_taylor_expansions}, which in our setting of a $C^{2}$ background has error $O(x^2)$, we get that:
  \begin{equation*}
    \abs{P_{11}+\kpp A}\lesssim \ep^{-\frac{1}{100}}e^{\kpp u}v,\qquad \abs{\p_xS_1-\kpp}\lesssim \ep^{-\frac{1}{100}}e^{\kpp u}v.
  \end{equation*}
  As a result, we obtain from the first two equations in \eqref{eq:double_null}
  \begin{align}
    \abs{\p_u \bp-A\kpp \bp}&\lesssim \ep^{-\frac{1}{100}}\abs{b_N}e^{\kpp (A+1) u}v^{N+1},\label{eq:b_bound}\\
    \abs{\p_v \ap}&\lesssim \abs{b_N}e^{\kpp(A+1)u}v^N\label{eq:a_bound}.
  \end{align}
  We will prove below analogous estimates, at top order, which will allow us to close the bootstrap estimates for $\ap,\bp)$.
  
  Applying the triangle inequality, we get that \eqref{eq:b_bound} implies
  \begin{equation*}
    \abs{\p_u \bp}\lesssim \abs{b_N}e^{\kpp A u}v^N.
  \end{equation*}
  We can then estimate
  \begin{align*}
    \left|\mrTht_1-\htTht_1+\kpp(\mrnu-\htnu)+\frac{\alp+1}{4(1-w_0)}\p_v \bp\right|&\lesssim \ep^{-\frac{2}{100}}\abs{b_N}e^{\kpp(A+1)u}v^N\\
    \left|\mrTht_2-\htTht_2\right|&\lesssim \ep^{-\frac{1}{100}}\abs{b_N}e^{\kpp A u}v^N.
  \end{align*}
  As a consequence, one obtains the following bounds
  \begin{align*}
    \left|\mrmu\mrTht_1-\htmu\htTht_1+\kpp(\mrnu-\htnu)+\frac{\alp+1}{4(1-w_0)}\p_v\bp\right|&\lesssim \ep^{-\frac{2}{100}}\abs{b_N} e^{\kpp(A+1)u}v^N,\\
    \left|\mrnu\mrTht_2-\htnu\htTht_2\right|&\lesssim \ep^{-\frac{1}{100}}\abs{b_N}e^{\kpp(A+1)u}v^N.
  \end{align*}
  Using the first bootstrap assumption in the first inequality above yields
  \begin{equation*}
    \left|\mrmu\mrTht_1-\htmu\htTht_1+\kpp(\mrnu-\htnu)+\frac{\alp+1}{4(1-w_0)}e^{\kpp A u}\p_v\brb\right|\lesssim \ep^{-\frac{2}{100}}\abs{b_N} e^{\kpp(A+1)u}v^N.
  \end{equation*}
  As a result, the third and fourth equations in \eqref{eq:double_null} imply that
  \begin{align}
    \left|\p_u (\mrnu-\htnu)-\kpp(\mrnu-\htnu)-\frac{\alp+1}{4(1-w_0)}e^{\kpp A u}\p_v\brb\right|&\lesssim \ep^{-\frac{2}{100}}\abs{b_N} e^{\kpp(A+1)u}v^N,\label{eq:nu_bound}\\
    \left|\p_v (\mrmu-\htmu)\right|&\lesssim \abs{b_N}e^{\kpp A u}v^{N-1}\label{eq:mubound}.
  \end{align}
  Integrating these equations gives
  \begin{align}
    \left|\mrnu-\htnu-\lmb(e^{\kpp Au}-e^{\kpp u})\p_v\brb\right|&\lesssim \ep^{-\frac{2}{100}}\abs{b_N}e^{\kpp (A+1)u}v^{N}\label{eq:nu_bootstrap},\\
    \abs{\mrmu-\htmu}&\lesssim \abs{b_N}e^{\kpp Au}v^{N}\label{eq:mu_bootstrap}.
  \end{align}
  In particular, we improved the third and fourth bootstrap assumption. We also note that using the bootstrap bounds and \eqref{eq:coordinate_evolution} we obtain
  \begin{equation*}
      \abs{\p_v \mrx-\p_v \htx}\lesssim \abs{\mrnu-\htnu}+e^{\kappa u}\left(\abs{\mrx-\htx}+|\ap|+|\bp|\right)\lesssim \abs{b_N}e^{\kpp A u}v^{N-1}.
  \end{equation*}
  By integration, we obtain that
  \begin{equation}
    \label{eq:x_bootstrap}
    \abs{\mrx-\htx}\lesssim \abs{b_N}e^{\kpp A u}v^N,
  \end{equation}
  so we improved the fifth bootstrap assumption as well.

  To conclude the bootstrap argument, we note that we can commute the $(\bp,\ap)$ transport equations with $(\p_v,\p_u)$ respectively
  , and obtain analogously to \eqref{eq:b_bound} and \eqref{eq:a_bound},
  \begin{align*}
    \abs{\p_u\p_v\bp-A\kpp \p_v \bp}&\lesssim \ep^{-\frac{1}{100}}\abs{b_N}e^{\kpp(A+1)u}v^N,\\
    \abs{\p_v\p_u\ap}&\lesssim \abs{b_N}e^{\kpp(A+1)u}v^N         .
  \end{align*}
    Integrating these equations, we obtain the following estimates  
  \begin{align}
    \abs{\p_v \bp-e^{\kpp A u}\p_v\brb}&\lesssim \ep^{-\frac{1}{100}}\abs{b_N}e^{\kpp (A+1) u}v^N\label{eq:b_bootstrap}\\
    \abs{\p_u \ap}&\lesssim \abs{b_N}e^{\kpp(A+1)u}v^{N+1}\label{eq:a_bootstrap}.
  \end{align}
  We note that the five estimates \eqref{eq:nu_bootstrap}, \eqref{eq:mu_bootstrap}, \eqref{eq:x_bootstrap}, \eqref{eq:b_bootstrap} and \eqref{eq:a_bootstrap} improve the given bootstrap assumptions, completing our bootstrap argument. At this point, we get by local existence theory that the solution of \eqref{eq:double_null} exists in $\calA$ and satisfies uniform $C^1$ bounds.

  \emph{Vanishing of $\mrnu$.}
  We now impose an additional smallness condition which is stronger than \eqref{eq:dlt_smallness_B}:
  \begin{equation*}
    \dlt L_{N+1}\leq \ep \abs{b_N}.
  \end{equation*}
  Using \eqref{eq:b_initialization}, we get that
  \begin{equation}
    \label{eq:b_der_bound}
    \left|\p_v\brb(\dlt)-b_N\frac{\dlt^{N-1}}{(N-1)!}\right|\lesssim \ep \abs{b_N}\dlt^{N-1}.
  \end{equation}
  From this, we obtain using \eqref{eq:nu_bootstrap} and \eqref{eq:upper_time}
  \begin{equation*}
    \left|\mrnu(u_+,\delta)-\htnu(u_+,\delta)+(1+\ep^{1-\frac{3}{100}})e^{\kpp u_+}\right|\lesssim \ep^{1-\frac{2}{100}} e^{\kpp u_+}.
  \end{equation*}
  In particular,
  \begin{equation*}
    \abs{e^{-\kpp u_+}\mrnu(u_+,\delta)+\ep^{1-\frac{3}{100}}}\lesssim \ep^{1-\frac{2}{100}}.
  \end{equation*}
  It follows that $e^{-\kpp u_+}\mrnu(u_+,\delta)<0$, and so $\mrnu(u_+,\delta)<0$. Thus, there must exist $u_S<u_{+}$ so that $\mrnu(u_S,\delta)=0$. One can similarly show that
  \begin{equation*}
    e^{-\kpp u_-}\mrnu(u_-,\dlt)-\ep^{1-\frac{3}{100}}\gtrsim \ep^{1-\frac{2}{100}},
  \end{equation*}
  which shows that $u_S>u_-$ and $\mrnu>0$ for $u\leq u_-$.
  
  We note that by \eqref{eq:b_der_bound}, we have $\p_v\bp<0$ whenever $v>0$. From the above estimates, one also obtains $\mrmu>0$ throughout the domain. We also note that one can show
  \begin{equation*}
      \kpp(A-1)+\p_u(e^{-\kpp u} \mrnu)|_{(u_S,\dlt)}\lesssim \ep^{1-\frac{3}{100}}
  \end{equation*}
  so $\partial_u\mrnu(u_S,\delta)<0,$ so our proof actually implies the existence of a $C^1$ curve where $\mrnu$ vanishes.
\end{proof}

We now prove Theorem \ref{thm:2}.
\begin{proof}
  By Lemma \ref{lem:null_shock_result}, the set $\calS_0$ defined in Remark \ref{rem:mghd} is non-empty and contains a point such that $\mrnu=0$. In addition, the quantities $\ap$ and $\bp$ are uniformly bounded. As $\mrL \bp=\mrnu^{-1}\p_v \bp$ and $\p_v \bp$ is bounded away from zero near this point, we get that $\mrL \bp$ blows up. We conclude that the solution of \eqref{eq:perturbation_eqn} is bounded in $L^{\infty}$ but it blows up in $C^{0,1}$.
\end{proof}

We comment that as $\dlt\to 0^+$, one has that $u_{\pm}\to \infty$. As such, the above proof yields a portion of the singular boundary of the solution which approaches the sonic cone as $v\to 0$ (see the meeting of green and blue curves in Figure \ref{fig:censorship_sketch}).

\subsection{Perturbations of Regularity $C^{\infty}$}
\label{sec:shock-form-c}

In our final theorem, we will consider background solutions and perturbation data which are $C^{\infty}$. In this case we will also drop the assumption that the solution is supported in the exterior. An appropriate definition of characteristic data in this case is given as follows.

\begin{definition}
\label{def:char_data_2}
    For the perturbation equations \eqref{eq:perturbation_eqn}, we define $C^{\infty}$ \emph{characteristic perturbation data} to consist of a point $p_-$ with $\tau(p_-)<0$, a $C^{\infty}$ curve $\calC_{p_-}$ emanating from $p_-$ passing through $O,$ and a function $\brb\in C^{\infty}(\calC_{p_-})$ vanishing to infinite order at $p_-$. We let $\ulcalC_{p_-}$ be a maximally extended $\underline{L}_b$ characteristic emanating from $p_-$ and impose $\ap|_{\ulcalC_{p_-}}=0$. We also impose that $\calC_{p_-}$ is tangent to $\mrL$.
\end{definition}

To continue the solution in an outer characteristic region up to the implosion time, we impose the following stability condition on the background.
\begin{definition}\label{def:admissible_global}
    Consider an admissible self-similar smooth implosion background as in Definition~\ref{def:admissible_background_smooth}, which is regular in the region $\left\{x\geq 0,\ t\leq0,\ u\geq 0\right\}\backslash\left\{(t,r)=(0,0)\right\}$. We say that it satisfies the \textit{exterior stability} condition if, for $X_0>0$ large enough, sufficiently small perturbations of the initial data on $\left\{x=X_0,\ t<0,\ u\geq 0\right\}\cup\left\{x\geq X_0,\ t<0,\ u=0\right\}$ yield solutions that remain regular in the region $\left\{x\geq X_0,\ t\leq0,\ u\geq 0\right\}\backslash\left\{(t,r)=(0,0)\right\}$.
\end{definition}
\begin{remark}
    We point out that we do not impose precise quantitative assumptions in the definition of the exterior stability condition since these depend on the background construction. In general, for self-similar smooth implosions, constructing the solution up to the implosion time already requires a notion of exterior stability as in Definition~\ref{def:admissible_global}. In particular, the smooth implosion examples in \cite{merleImplosionCompressibleFluid2022e}, \cite{buckmasterSmoothImplodingSolutions2025b} and \cite{shao2025blowup3dcompressiblenavierstokes} satisfy the exterior stability condition.
\end{remark}

We now state our third main result.
\begin{theorem}
  \label{thm:3}
  
  \begin{enumerate}
    \item 
      Consider an admissible self-similar smooth implosion background as in Definition~\ref{def:admissible_background_smooth}. Consider the perturbation equations \eqref{eq:perturbation_eqn} with $C^{\infty}$ perturbation data, as in Definition \ref{def:char_data_2}. Suppose that the perturbation data satisfies
  \begin{equation}
      \label{eq:derivative_conditions}
      \mrL^{(k)}\brb|_{p_-}=0,\qquad\text{ for all }k\geq0.
  \end{equation}
    Let $\dlt=-\tau(p_-)$ and let $p_+$ the point along $\calC_{p_-}$ with $\tau(p_+)=\dlt$. Let $\calC_{p_-}'\subset \calC_{p_-}$ denote the subcurve from $p_-$ to $p_+$. We fix an integer $N\geq 1$, and define
    \begin{equation*}
        b_1:=\inf_{\calC_{p_-}'} \mrL \brb,\qquad B_N:=\left\|(\mrL)^{(N)}\brb\right\|_{L^{\infty}(\calC_{p_-}')},\qquad \Lmb:=\frac{(2\dlt)^{N-1}B_N}{\abs{b_1}(N-1)!}.
    \end{equation*}
    Suppose that 
    \begin{equation}
        \label{eq:b1_bound}
        0<-b_1\leq \lmb^{-1}
    \end{equation}
    and for some implicit constant depending on the background
    \begin{equation}\label{eq:smallness_implicit}
        \frac{\dlt \Lmb^2}{\abs{b_1}^{\frac{1}{A-1}}}\ll 1.
    \end{equation}
    Then, in the domain of dependence of $\calC_{p_-}',$ the solution blows up in $C^{0,1}$ while remaining bounded in $L^{\infty}$.
\item Suppose that we additionally have
    \begin{equation}
        \label{eq:sign_ineq}
        \mrL \brb|_{p_+}= 0.
    \end{equation}
   Let $\ulcalC_{p_+}$ denote the ingoing characteristic from $p_+$, contained in the MGHD of the characteristic data. Then, a shock forms on the boundary of the MGHD, away from $\ulcalC_{p_+}$. The solution remains $C^1$ regular along $\ulcalC_{p_+}$ and satisfies the following bounds along $\ulcalC_{p_+}$:
\begin{equation*}
    \abs{\ap}\lesssim \Lambda B_N e^{\kpp(A+1)u}\dlt^{N+1},\qquad \abs{\mrulL \ap}\lesssim \Lambda B_N e^{\kpp(A+1)u}\dlt^{N+1},\qquad x\lesssim \Lambda e^{\kpp u}\dlt.
\end{equation*}
\begin{equation*}
    \abs{\bp}\lesssim B_N e^{\kpp A u}\dlt^N,\qquad \abs{\mrL \bp}\lesssim B_N e^{\kpp(A-1)u}\dlt^{N-1}.
\end{equation*}
  \item Suppose that in addition to the conditions in (a) and (b), the background satisfies the exterior stability condition in Definition~\ref{def:admissible_global}. We extend the outgoing characteristic $\calC_{p_-}$ up to $\{t=0\}$, and we denote by $\calC_{p_+}\subset \calC_{p_-}$ the subcurve to the future of $p_+$. We suppose that the $C^{\infty}$ perturbation data satisfy
  \begin{equation}
      \brb=0,\qquad\text{ for all }\tau\geq\delta.
  \end{equation}  
  Then the solution of the perturbation equations \eqref{eq:perturbation_eqn} exists and remains smooth in the domain of dependence of $\ulcalC_{p_+}\cup\calC_{p_+}$, up to the time of implosion.
  \end{enumerate}
\end{theorem}
  We comment that the first condition in \eqref{eq:derivative_conditions} ensures that the data may be smoothly glued to the interior. One can choose $\brb$ to have compact support, which implies \eqref{eq:sign_ineq}.

  We note the following corollary (cf. Corollary \ref{cor:instability_corollary} with $k=\lfloor A\rfloor-1$).
\begin{corollary}{(Short pulse instabilities).}
\label{cor:smooth_pulse}
    Consider a nontrivial fixed profile $\chi\in C_0^{\infty}$ for $\brb$. There exists a sequence $(t_j,\dlt_j)\to 0$ so that the rescaled profiles $t_j\cdot \chi(\cdot/\dlt_j)$ satisfy the hypotheses of Theorem \ref{thm:3} and converge to $0$ in $C^{\lfloor A\rfloor}$.
\end{corollary}
\begin{proof}
    We note that $\Lmb$ is constant for the sequence, so \eqref{eq:smallness_implicit} becomes
    \begin{equation*}
       \dlt^A\ll t.
   \end{equation*}
    This is satisfied if one picks $\dlt=(\ep t)^{\frac{1}{A}},$ for some fixed $\ep>0$ depending on the profile. For every integer $k\geq0$, the $\dot{C}^{k}$ norm is of size $t\cdot \dlt^{-k}=t^{1-\frac{k}{A}}\ep^{-\frac{k}{A}}$, and so by taking $t\to 0$, we see that if $k<A$, the $\dot{C}^k$ norm goes to $0$. We note that as $\chi$ is compactly supported, there must exist a point where $\mrL \chi<0$. In addition, as $t\to 0$, we have that $-b_1\sim t^{1-\frac{1}{A}}\ep^{-\frac{1}{A}}\to 0$, so \eqref{eq:b1_bound} is satisfied. 
\end{proof}
\begin{remark}
\label{rem:C_infinity_smallness}
  The previous corollary is sharp in the sense that the data for Theorem \ref{thm:3} cannot be chosen to be small in the $\dot{C}^{\lfloor A\rfloor+1}$ topology. To see this, note that by integrating $\p_v^{N}b$ from $p_-$ a total of $N-1$ times, we see $\Lmb\gtrsim 1$. In particular, as we consider a bounded range of $\abs{b_1}$, it is necessary for \eqref{eq:smallness_implicit} that $\dlt\ll 1$. Similarly, since
  \begin{equation*}
      \frac{1}{\dlt^{N-1}B_N}\lesssim \frac{1}{\abs{b_1}},
  \end{equation*}
  we note that a necessary condition for \eqref{eq:smallness_implicit} is
  \begin{equation*}
    \dlt^{A-N}\ll B_N.
  \end{equation*}
  In particular, if $N>A$ a necessary condition for this statement is that $B_{N}\gg 1$.
\end{remark}

Similarly to the proof of Theorem \ref{thm:2}, we will prove this result at the end of the section as a consequence of the following result for the system of null perturbation equations. For this system, we will now consider characteristic curves which intersect the interior. We will have that $p_-$ corresponds to $(0,-\dlt)$ in the $(u,v)$ coordinates. We define
\begin{equation}
    \label{eq:ulcalC_def_null_2}
    \ulcalC_{p_-}:=\left\{(u,v)\colon v=-\delta,0\leq u<u_0\right\},\qquad \calC_{p_-}:=\left\{(u,v)\colon u=0,-\dlt \leq v\leq \Delta\right\},
\end{equation}
    where the parameter $\dlt$ satisfies $0<\dlt<\Dlt$, and $u_0$ is defined so that $\lim_{u\rightarrow u_0}\htx(u,-\delta)=-\infty$.

\begin{lemma}\label{lem:null_shock_smooth}
    \begin{enumerate}
    \item 
        Consider an admissible self-similar smooth implosion background as in Definition~\ref{def:admissible_background_smooth}. We define $\ulcalC_{p_-}$, $\calC_{p_-}$ by \eqref{eq:ulcalC_def_null_2} with $\dlt\in (0,\Dlt)$, and we consider characteristic data $\ap|_{v=-\dlt}=0$ and $\bp|_{u=0}=\brb(v)$ for the null perturbation system \eqref{eq:double_null}. We assume that $\brb\in C^{\infty}(\calC_{p_-})$ satisfies for all $k\in \N$
    \begin{equation*}
        \p_v^{k}\brb|_{v=-\dlt}=0.
    \end{equation*}
    We impose the gauge conditions in Subsection~\ref{sec:pertubation_gauge} with $\delta'=\delta$:
    \begin{equation*}
    \mrmu|_{\ulcalC_{p_-}}=\htmu|_{\ulcalC_{p_-}},\qquad \mrnu|_{\calC_{p_-}}=1,\qquad \mrx|_{\ulcalC_{p_-}}=\htx|_{\ulcalC_{p_-}}.
    \end{equation*}
    Let $\calC_{p_-}'\subset \calC_{p_-}$ be a subcurve from $p_-$ to $p_+$, which corresponds to $(0,\dlt)$ in the $(u,v)$ coordinates. We fix an integer $N\geq1$, and define
    \begin{equation*}
        b_1:=\inf_{\calC_{p_-}'} \p_v\brb,\qquad B_N:=\left\|\p_v^N\brb\right\|_{L^{\infty}(\calC_{p_-}')},\qquad \Lmb:=\frac{(2\dlt)^{N-1}B_N}{\abs{b_1}(N-1)!}.
    \end{equation*}
    Suppose that $0<-\lmb b_1\leq 1$ and for some implicit constant depending on the background
    \begin{equation}\label{eq:smallness_implicit_lemma}
        \frac{\dlt \Lmb^{2}}{\abs{b_1}^{\frac{1}{A-1}}}\ll 1.
    \end{equation}
    In this case, there exists explicit $u_-,u_+$ depending on the data so that if we define
  \begin{equation*}
    \calA:=[0,u_+)\times [-\dlt,\dlt],
  \end{equation*}
  then the solution of \eqref{eq:double_null} exists in $\calA$ and satisfies $\mrmu>0$. The level set $\{\mrnu=0\}$ is a subset of $\{\p_v\bp<0\}$ and for any $v_*$ such that $\p_v \brb|_{v_*}=b_1$ this level set contains $(u_S,v_*)$ for some $u_S\in (u_-,u_+)$.
  \item Suppose that in addition to the conditions in (a), we also have that
  \begin{equation}
  \label{eq:derivative_vanishing}
      \p_v\brb|_{v=\dlt}=0.
  \end{equation}
  Then the minimum $u$ value along $\{\mrnu=0\}$ is achieved at a point $(u_M,v_M)\in(u_-,u_+)\times(-\dlt,\dlt)$, where $u_{\pm}$ is defined in \eqref{eq:def_upm}. We denote the ingoing characteristic passing through $p_+$ by
  \begin{equation*}
      \ulcalC_{p_+}:=\left\{(u,v)\colon 0\leq u<u_M,v=\dlt\right\}.
  \end{equation*}
  Along $\ulcalC_{p_+}$, we have that the solution satisfies the following estimates
  \begin{equation*}
      \abs{e^{-\kpp u}\mrnu-1}\ll 1,\qquad \abs{\bp}\lesssim B_N e^{\kpp A u}\dlt^{N},\qquad \abs{\p_v \bp}\lesssim B_N e^{\kpp A u}\dlt^{N-1},
  \end{equation*}
  \begin{equation*}
  \abs{\mrmu-1}\ll 1,\qquad
      \abs{\ap}\lesssim \Lambda B_N e^{\kpp(A+1)u}\dlt^{N+1},\qquad \abs{\p_u \ap}\lesssim \Lambda B_N e^{\kpp(A+1)u}\dlt^{N+1}.
  \end{equation*}
  \item  In addition to the conditions in (a) and (b), we assume that the background satisfies the exterior stability condition in Definition~\ref{def:admissible_global}. Moreover, we also assume that $\brb$ vanishes to infinite order at $v=\delta.$ We extend the outgoing characteristic $\calC_{p_-}$ for all $v\geq\delta$, and we extend the characteristic data and gauge conditions on $\calC_{p_-}$ for $v\geq\delta$ as
  \[\bp|_{u=0}(v)=0,\qquad \mrnu|_{u=0}(v)=1.\]
  Then the solution of \eqref{eq:double_null} with the given characteristic data on $\calC_{p_-}$ and the induced characteristic data on $\ulcalC_{p_+}$ exists and remains smooth in the region $\left\{0\leq u< u_M,\ \delta\leq v,\ t<0\right\},$ with $\mrnu, \mrmu>0.$
  \end{enumerate}
\end{lemma}
\begin{proof}
  The proof of part (a) follows a similar strategy to that of Lemma \ref{lem:null_shock_result}. For part (b), the previous estimates allow us to conclude that, under the assumptions \eqref{eq:derivative_vanishing}, $e^{-\kpp u}\mrnu$ must be close to $1$ along $\ulcalC_{p_+}$, which yields the desired statements. Finally, in part (c) we use the exterior stability property of the background solution to conclude existence and smoothness in the outer characteristic region, up to the time of the implosion. A sketch of the characteristic domain, under the assumption of (b), is included in Figure \ref{fig:characteristic_domain_2}. Note that in this result, the singular boundary of the MGHD does not meet $\ulcalC_{p_+}$ (cf. Figure \ref{fig:characteristic_domain}).

    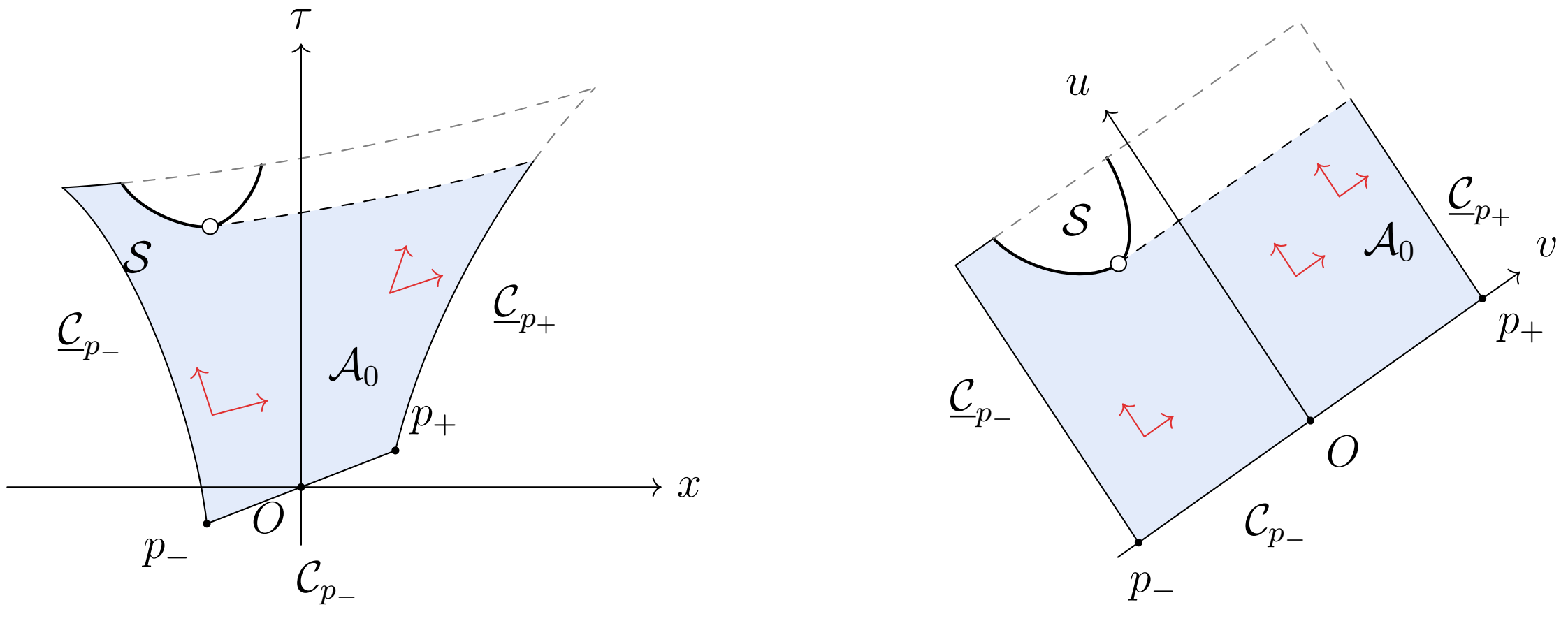
\begin{figure}[h]
    \centering
\begin{tikzpicture}[scale=1.05,line width=0.4pt,
    line cap=round,line join=round]
    \definecolor{domainblue}{RGB}{40,103,224}
    \definecolor{domainred}{RGB}{225,48,48}
    \tikzset{
        singular/.style={black,line width=0.8pt},
        continuation/.style={dashed,black!50},
        horizon/.style={dashed,black},
        causal/.style={->,domainred}
    }

    \begin{scope}
        \coordinate (LO) at (0,0);
        \coordinate (Lm) at (-0.85,-0.33);
        \coordinate (Lp) at (0.85,0.33);
        \coordinate (LT) at (-2.15,2.70);
        \coordinate (LF) at (2.65,3.60);
        \coordinate (LB) at (-1.617718,2.741852);
        \coordinate (LM) at (-0.82,2.35);
        \coordinate (LE) at (-0.357221,2.903906);
        \coordinate (LH) at (2.100277,2.945190);

        \fill[domainblue,fill opacity=0.13]
            (Lm) -- (LO) -- (Lp)
            .. controls (1.087600,1.316400)
                    and (1.605136,2.266512) .. (LH)
            .. controls (1.126851,2.653162)
                    and (0.153426,2.486280) .. (LM)
            .. controls (-1.02,2.322)
                    and (-1.457718,2.481852) .. (LB)
            .. controls (-1.798400,2.724768)
                    and (-1.976000,2.710800) .. (LT)
            .. controls (-1.55,2.20) and (-0.97,0.72) .. (Lm)
            -- cycle;

        \draw (Lm) -- (LO) -- (Lp);

        \draw (Lm) .. controls (-0.97,0.72)
                           and (-1.55,2.20) .. (LT);
        \draw (Lp) .. controls (1.087600,1.316400)
                           and (1.605136,2.266512) .. (LH);

        \draw (LT) .. controls (-1.976000,2.710800)
                           and (-1.798400,2.724768) .. (LB);

        \draw[singular]
            (LB) .. controls (-1.457718,2.481852)
                        and (-1.02,2.322) .. (LM)
            .. controls (-0.64,2.3752)
                    and (-0.407221,2.62) .. (LE);

        \draw[horizon]
            (LM) .. controls (0.153426,2.486280)
                        and (1.126851,2.653162) .. (LH);

        \draw[continuation]
            (LB) .. controls (-0.292720,2.867136)
                        and (1.198000,3.160000) .. (LF);
        \draw[continuation]
            (LH) .. controls (2.292832,3.209120)
                        and (2.482000,3.432000) .. (LF);

        \draw[causal] (-0.80,0.65) -- (-0.94,1.08);
        \draw[causal] (-0.80,0.65) -- (-0.30,0.78);
        \draw[causal] (0.80,1.75) -- (0.95,2.18);
        \draw[causal] (0.80,1.75) -- (1.28,1.91);

        \draw[->] (-2.65,0) -- (3.25,0) node[right] {$x$};
        \draw[->] (0,-0.52) -- (0,4.00) node[above] {$\tau$};
        \foreach \point in {Lm,LO,Lp}
            \fill (\point) circle (1pt);
        \filldraw[fill=white] (LM) circle (2pt);

        \node[below left] at (Lm) {$p_-$};
        \node[below left] at (LO) {$O$};
        \node[above right] at (Lp) {$p_+$};
        \node at (0.25,-0.87) {$\mathcal C_{p_-}$};
        \node[anchor=east] at (-1.45,1.35)
            {$\underline{\mathcal C}_{p_-}$};
        \node[anchor=west] at (1.60,1.60)
            {$\underline{\mathcal C}_{p_+}$};
        \node at (0.48,1.08) {$\mathcal A_0$};
        \node[anchor=north east] at (-1.20,2.35) {$\mathcal S$};
    \end{scope}

    \begin{scope}[shift={(9.10,0.60)},
        x={(-1.65cm,2.50cm)},y={(1.55cm,1.10cm)}]

        \fill[domainblue,fill opacity=0.13]
            (0,-1) -- (0,1) -- (0.72,1) -- (0.72,-0.35)
            .. controls (0.72,-0.52) and (0.86,-0.74) .. (1,-0.78)
            -- (1,-1) -- cycle;

        \draw (0,-1) -- (1,-1) -- (1,-0.78);
        \draw (0,1) -- (0.72,1);

        \draw[singular]
            (1,-0.78) .. controls (0.86,-0.74)
                            and (0.72,-0.52) .. (0.72,-0.35)
            .. controls (0.72,-0.22) and (0.90,-0.12) .. (1,-0.12);

        \draw[horizon] (0.72,-0.35) -- (0.72,1);
        \draw[continuation] (1,-0.78) -- (1,1) -- (0.72,1);

        \foreach \u/\v in {0.25/-0.70,0.38/0.32,0.50/0.70} {
            \draw[causal] (\u,\v) -- ({\u+0.12},\v);
            \draw[causal] (\u,\v) -- (\u,{\v+0.17});
        }

        \draw[->] (0,0) -- (1.12,0) node[above left] {$u$};
        \draw[->] (0,-1.12) -- (0,1.22) node[above right] {$v$};
        \foreach \v in {-1,0,1}
            \fill (0,\v) circle (1pt);
        \filldraw[fill=white] (0.72,-0.35) circle (2pt);

        \node[below,xshift=4pt,yshift=-4pt] at (0,-1) {$p_-$};
        \node[below right] at (0,0) {$O$};
        \node[below right] at (0,1) {$p_+$};
        \node[anchor=north west,yshift=-3pt] at (0,-0.47)
            {$\mathcal C_{p_-}$};
        \node[anchor=east,xshift=-4pt] at (0.50,-1)
            {$\underline{\mathcal C}_{p_-}$};
        \node[anchor=west,xshift=4pt] at (0.35,1)
            {$\underline{\mathcal C}_{p_+}$};
        \node at (0.30,0.78) {$\mathcal A_0$};
        \node at (0.90,-0.40) {$\mathcal S$};
    \end{scope}
\end{tikzpicture}
    \caption{Sketch of the characteristic domain $\mathcal{A}$ for Theorem~\ref{thm:3}.}
    \label{fig:characteristic_domain_2}
\end{figure}

  As in the proof of Lemma~\ref{lem:null_shock_result}, we use $``\lesssim"$ to denote inequalities with implicit constants which only depend on the background, and we consider $\ep>0$ to be a sufficiently small quantity depending on the background.
  We begin by defining the timescales $u_{\pm}$ by
  \begin{equation}\label{eq:def_upm}
    \lmb b_1(e^{\kpp(A-1)u_\pm}-1)=-1\mp \ep^{1-\frac{3}{100}}.
  \end{equation}
  Our hypothesis is that the following quantitative assumptions hold
  \begin{equation}\label{eq:first_estimate}
    \frac{\dlt \Lmb^2}{\abs{b_1}^{\frac{1}{A-1}}}\leq\ep,
  \end{equation}
  \begin{equation}\label{eq:second_estimate}
    \lmb \abs{b_1}\leq 1.
  \end{equation}
  We comment that the estimate \eqref{eq:second_estimate} will ensure that shock formation is driven by the blue-shift amplification. We note that \eqref{eq:second_estimate} implies
  \begin{equation}\label{eq:third_estimate}
    \lmb \abs{b_1}e^{\kpp(A-1)u_+}\leq 3.
  \end{equation}
  As a consequence of \eqref{eq:third_estimate} and the fact that $\Lambda\geq1$, we obtain the following smallness estimate (in fact, we obtain the stronger estimate \eqref{eq:improved_smallness}, that we do not need until proving shock formation)
  \begin{equation}
  \label{eq:thm_3_smallness}
    \dlt e^{\kpp u_+}\Lmb\lesssim \ep.
  \end{equation}
    Moreover, the definition of $\Lambda$ implies that we have the initial data estimate
  \begin{equation*}
    \norm{\p_v \brb}_{L^{\infty}(\mathcal{C}')}\leq B_N\frac{\delta^{N-1}}{(N-1)!}=\Lmb \abs{b_1}.
  \end{equation*}

  We note that in particular, this implies that $\abs{\mrx}\lesssim \ep^{1-\frac{1}{100}}$. We comment that $r=0$ corresponds to $\mrx=-\infty$, and so these hypotheses ensure that the domain we consider is far from the axis of symmetry. Similarly, the asymptotic expansions in Subsection \ref{sec:append-null-coord} ensure that $\abs{e^{-\kpp u}\htnu(u,-\dlt)-1}\lesssim \ep$ and hence the conversion to null coordinates is well-behaved along $v=-\dlt$.

    Similarly to the proof of Lemma \ref{lem:null_shock_result}, we carry out a bootstrap argument. We consider the following bootstrap assumptions in this case
    \begin{align*}
        \abs{\p_v \bp-e^{\kpp A u}\p_v \brb}&\leq \ep^{-\frac{2}{100}}\Lmb^2\abs{b_1}e^{\kpp (A+1) u}\dlt,\\
        \abs{\p_u\ap}&\leq \ep^{-\frac{1}{100}}\Lmb^2 \abs{b_1}e^{\kpp(A+1)u}\dlt^2,\\
        \abs{\mrnu-\htnu-\lmb(e^{\kpp A u}-e^{\kpp u})\p_v \brb}&\leq \ep^{-\frac{3}{100}}\Lmb^2\abs{b_1}e^{\kpp(A+1)u}\dlt,\\
        \abs{\mrmu-\htmu}&\leq \ep^{-\frac{1}{100}}\Lmb \abs{b_1}e^{\kpp A u}\dlt,\\
        \abs{\mrx-\htx}&\leq \ep^{-\frac{1}{100}}\Lmb \abs{b_1}e^{\kpp A u}\dlt.
    \end{align*}
    Once again, one can check that the bootstrap assumptions hold initially with significantly better constants. We note that the bootstrap assumptions imply
    \begin{equation*}
    \abs{\bp-e^{\kpp A u}\brb}\lesssim \ep^{-\frac{2}{100}}\Lmb^2 \abs{b_1}e^{\kpp(A+1)u}\dlt^2,\qquad \abs{\ap}\lesssim \ep^{-\frac{1}{100}}\Lmb^2 \abs{b_1}e^{\kpp(A+1)u}\dlt^2,
    \end{equation*}
    \begin{equation*}
    \abs{\p_v \bp}\lesssim \Lmb \abs{b_1}e^{\kpp A u},\qquad \abs{\bp}\lesssim \Lmb \abs{b_1}e^{\kpp A u}\dlt,\qquad \abs{\mrnu-\htnu}\lesssim \Lmb \abs{b_1}e^{\kpp A u}.
    \end{equation*}
  As before, we will use below the following nonlinear bounds
  \begin{equation*}
    \abs{\bp^2}\lesssim \Lmb^2 \abs{b_1}e^{\kpp(A+1)u}\dlt^2,\  \abs{\bp\p_v \bp}\lesssim \Lmb^2 \abs{b_1}e^{\kpp(A+1)u}\dlt,\  \abs{(\mrnu-\htnu)\bp}\lesssim \Lmb^2 \abs{b_1}e^{\kpp(A+1)u}\dlt.
  \end{equation*}
  \emph{Closing the bootstrap argument.}
  From the first two equations in \eqref{eq:double_null}, we obtain
  \begin{align}
    \abs{\p_u \bp-A\kpp \bp}&\lesssim \ep^{-\frac{1}{100}}\Lmb^2 \abs{b_1}e^{\kpp(A+1)u}\dlt^2,\label{eq:b_bound_3}\\
    \abs{\p_v \ap}&\lesssim \Lmb^2 \abs{b_1}e^{\kpp (A+1)u}\dlt.\label{eq:a_bound_3}
  \end{align}
  As a consequence of \eqref{eq:b_bound_3}, we get
  \begin{equation*}
    \abs{\p_u \bp}\lesssim \Lmb \abs{b_1}e^{\kpp A u}\dlt.
  \end{equation*}
  Similarly to the derivation of \eqref{eq:nu_bound}-\eqref{eq:mubound}, we obtain the estimates
  \begin{align}
    \left|\p_u(\mrnu-\htnu)-\kpp(\mrnu-\htnu)-\frac{\alp+1}{4(1-w_0)}e^{\kpp A u}\p_v\brb\right|&\lesssim \ep^{-\frac{2}{100}}\Lmb^2 \abs{b_1}e^{\kpp(A+1)u}\dlt,\\
    \left|\p_v(\mrmu-\htmu)\right|&\lesssim \Lmb \abs{b_1}e^{\kpp A u}.
  \end{align}
  We integrate these equations to get
  \begin{align}
    \left|\mrnu-\htnu-\lmb(e^{\kpp A u}-e^{\kpp u})\p_v \brb\right|&\lesssim \ep^{-\frac{2}{100}}\Lmb^2 \abs{b_1}e^{\kpp(A+1)u}\dlt,\label{eq:nu_bootstrap_3}\\
    \left|\mrmu-\htmu\right|&\lesssim \Lmb \abs{b_1}e^{\kpp A u}\dlt.\label{eq:mu_bootstrap_3}
  \end{align}
  In particular, we improved the third and fourth bootstrap assumption. Using \eqref{eq:coordinate_evolution}, we also obtain the bound
  \begin{equation*}
    \abs{\p_v \mrx-\p_v\htx}\lesssim \Lmb \abs{b_1}e^{\kpp A u},
  \end{equation*}
  which upon integration yields
  \begin{equation}
    \abs{\mrx-\htx}\lesssim \Lmb \abs{b_1}e^{\kpp A u}\dlt\label{eq:x_bootstrap_3},
  \end{equation}
  improving the fifth bootstrap assumption as well.

  To conclude the bootstrap argument, we differentiate the transport equations for $\bp$ and $\ap$, and compute the following estimates
  \begin{align*}
    \abs{\p_u\p_v \bp-A \kpp \p_v \bp}\lesssim \ep^{-\frac{1}{100}}\Lmb^2 \abs{b_1}e^{\kpp(A+1)u}\dlt\\
    \abs{\p_v\p_u \ap}\lesssim \Lmb^2 \abs{b_1}e^{\kpp(A+1)u}\dlt.
  \end{align*}
  Integrating these equations, we obtain
  \begin{align}
    \abs{\p_v \bp-e^{\kpp A u}\p_v\brb}&\lesssim \ep^{-\frac{1}{100}}\Lmb^2 \abs{b_1}e^{\kpp(A+1)u}\dlt,\label{eq:b_bootstrap_3}\\
    \abs{\p_u \ap}&\lesssim \Lmb^2 \abs{b_1}e^{\kpp(A+1)u}\dlt^2.\label{eq:a_bootstrap_3}
  \end{align}
  The five bounds \eqref{eq:nu_bootstrap_3}, \eqref{eq:mu_bootstrap_3}, \eqref{eq:x_bootstrap_3}, \eqref{eq:b_bootstrap_3}, and \eqref{eq:a_bootstrap_3} improve the given bootstrap assumptions, completing our bootstrap argument. At this point, we get by local existence theory that the solution of \eqref{eq:double_null} exists in $\calA$ and satisfies uniform $C^1$ bounds.
  
  \emph{Vanishing of $\mrnu$.} We recall that \eqref{eq:first_estimate} and \eqref{eq:third_estimate} actually imply the following smallness estimate (which is stronger than \eqref{eq:thm_3_smallness})
  \begin{equation}
  \label{eq:improved_smallness}
    \dlt e^{\kpp u_+}\Lmb^2\lesssim \ep.
  \end{equation}
  It follows that along $\{u=u_+\}$, \eqref{eq:nu_bootstrap_3} and the definition of $u_+$ imply for all $v\in[-\dlt,\delta]$
  \begin{equation*}
    \left|\mrnu(u_+,v)-\htnu(u_+,v)+(1+\ep^{1-\frac{3}{100}})e^{\kpp u_+}\frac{\p_v \brb(v)}{b_1}\right|\lesssim \ep^{1-\frac{2}{100}} e^{\kpp u_+}.
  \end{equation*}
  Thus, we get that for all $v\in[-\dlt,\delta]$
  \begin{equation}
    \label{eq:mrnu_estimate}
    \left|e^{-\kpp u_+}\mrnu(u_+,v)+(1+\ep^{1-\frac{3}{100}})\frac{\p_v \brb(v)}{b_1}-1\right|\lesssim \ep^{1-\frac{2}{100}}.
  \end{equation}
  By assumption, there exists $v_*\in(-\dlt,\delta]$ such that $\p_v\brb(v_*)=b_1$, so we have $\mrnu(u_+,v_*)<0$. Thus, there exists $u_S<u_+$ so that $\mrnu(u_S,u_*)=0$. We also note that it follows from a similar estimate, using $\p_v\brb\geq b_1$, that
  \begin{equation*}
    \mrnu(u_-,v)-\ep^{1-\frac{3}{100}}\gtrsim \ep^{1-\frac{2}{100}},
  \end{equation*}
  so that $\mrnu$ is non-vanishing for $u=u_-$ and in particular $u_-<u_S$.
  We then note that by a slight variation of \eqref{eq:nu_bootstrap_3}, replacing $\partial_v\brb$ by $e^{-\kappa A u}\partial_vb$, we have that for all $(u,v)\in\mathcal{A}$
  \begin{equation*}
    e^{-\kpp u}\mrnu\geq 1-\ep^{1-\frac{3}{100}}+\lambda\left(e^{-\kappa u}-e^{-\kappa A u}\right)\partial_v\bp.
  \end{equation*}
  In particular, whenever $\p_v \bp\geq 0$, we have that
  \begin{equation*}
    e^{-\kpp u}\mrnu\geq 1-\ep^{1-\frac{3}{100}}> 0,
  \end{equation*}
  which proves the last claim in part (a).

  \emph{Part (b).}
  In this case, we note that since $\p_v \brb|_{\dlt}= 0,$ one can similarly conclude from \eqref{eq:nu_bootstrap_3} that at along $\{v=\dlt\},$ we have for all $u\in[0,u_+]$
  \begin{equation}
      \label{eq:nu_smallness}
      \left|e^{-\kpp u}\mrnu-1\right| \lesssim \ep^{1-\frac{2}{100}}.
  \end{equation}
  It follows that $\{\mrnu=0\}$ does not intersect $\{v=\dlt\},$ which implies that the minimum value of $u$ occurs at a point $(u_M,v_M)$, for some $u_M\in(u_-,u_+)$ and $v_M\in (0,\dlt)$. The stated bounds follow from \eqref{eq:nu_smallness} and the previous estimates.

  \emph{Part (c).} We first show the solution satisfies improved estimates on $\left\{0\leq u< u_M,\ v=\delta\right\},$ due to the vanishing of $\overline{b}_p$ and $\partial_v\overline{b}_p$ at $v=\delta.$ Thus, we get that:
    \begin{align}
        \abs{a_p}+\abs{\p_u a_p}&\lesssim\Lmb^2 \abs{b_1}e^{\kpp(A+1)u}\dlt^2,\qquad \abs{\p_v a_p}\lesssim\Lmb^2 \abs{b_1}e^{\kpp(A+2)u}\dlt^2,\label{eq:ext1}\\ 
        \abs{b_p}+\abs{\p_u b_p}&\lesssim\Lmb^2 \abs{b_1}e^{\kpp(A+1)u}\dlt^2,\qquad \abs{\p_v b_p}\lesssim\Lmb^2 \abs{b_1}e^{\kpp(A+2)u}\dlt^2,\label{eq:ext2}\\
        \abs{\mrx-\htx}&\lesssim \Lmb^2 \abs{b_1}e^{\kpp(A+1)u}\dlt^2,\qquad
        \abs{e^{-\kpp u}\mrnu-1}+\abs{\mrmu-1}\lesssim \dlt e^{\kpp u}.\label{eq:ext3}
    \end{align}

    Next, we notice that $\mrx$ is lower bounded up to a constant by $\delta$ on $\underline{\mathcal{C}}_{p_+}$. We consider $x_0>0$ small, depending on the background (but independent of $\delta$ and $\ep$), so that the Taylor expansions in Section~\ref{sec:2} hold for $x\in[0,x_0].$ A similar argument to the proof of Lemma~\ref{lem:null_shock_smooth}, together with the derivation of \eqref{eq:ext1}-\eqref{eq:ext3} using the vanishing of $\bp|_{u=0}$ for $v\geq\delta,$ shows that the solution exists and satisfies $C^1$ estimates in the region $\left\{0\leq u< u_M,\ \delta\leq v\right\}\cap\left\{\mrx\leq x_0\right\}.$ In particular, in terms of the $(\tau,x)$ coordinates, the following bound holds on $\left\{x=x_0,\ 0\leq u<u_M\right\}$
    \begin{equation}\label{eq:bounds_x_0}
        \abs{a_p}+\abs{\p_x a_p}+\abs{\p_{\tau} a_p}+\abs{b_p}+\abs{\p_x b_p}+\abs{\p_{\tau} b_p}\lesssim\ep^2.
    \end{equation}
    Moreover, commuting the equations yields higher order analogues of \eqref{eq:bounds_x_0}.

    Since the background solution satisfies the exterior stability condition, we consider $X_0>0$, depending on the background, as in Definition~\ref{def:admissible_global}. The region $\left\{0\leq u< u_M,\ \delta\leq v\right\}\cap\left\{x_0\leq\mrx\leq X_0\right\}$ is finite in the self-similar variable $x,$ so it is straightforward to show that \eqref{eq:bounds_x_0} also holds on $\left\{x=X_0,\ 0\leq u<u_M\right\}$, with worse constants.
    
    Finally, starting from the induced data on the hypersurface $\left\{x=X_0,\ 0\leq u< u_M\right\}\cup\left\{u=0,\ x\geq X_0\right\}$, the exterior stability condition implies that the solution exists and is smooth in the remaining region $\left\{0\leq u< u_M,\ x\geq X_0,\ t< 0\right\}.$
\end{proof}

We now prove Theorem \ref{thm:3}.
\begin{proof}
      Part (a) of this result follows from Lemma~\ref{lem:null_shock_smooth} (similarly to the proof of Theorem~\ref{thm:2} and Lemma~\ref{lem:null_shock_result}). For part (b), we note that the condition $\abs{e^{-\kpp u}\mrnu-1}\ll 1$ on $\{v=\dlt\}$ allows one to convert from the $(u,v)$ coordinates back to the $(\tau,x)$ coordinates using \eqref{eq:transport_renormalizations}. Similarly, for part (c) we use the fact that $\mrmu,\mrnu>0,$ to obtain the desired result in the $(\tau,x)$ coordinates. Finally, the exterior stability condition implies that the solution extends smoothly to $\{t=0,\ 0\leq u<u_M\}$, with respect to the Cartesian radial coordinates $(t,r)$.
\end{proof}

\printbibliography
\end{document}